\documentclass[11pt]{article}
\usepackage{amsmath,amssymb}
\usepackage{graphicx}
\usepackage{amscd}
\usepackage{amsthm}

\makeatletter
\@addtoreset{equation}{section}
\makeatother

\newcommand{\Z}{{\mathbb{Z}}}

\def \tei{\widetilde{e_i}}
\def \tfi{\widetilde{f_i}}
\def \te0{\widetilde{e_0}}
\def \tf0{\widetilde{f_0}}
\def \gge{\mathfrak{g}}
\def \uge{U_q(\mathfrak{g})}
\def \udn{U_q(D(N))}
\def \ovB{\overline{B}}
\def \ovV{\overline{V}}
\def \lam{\lambda}
\def \Lam{\Lambda}
\def \omg{\omega}
\def \ome0{\omega_0}
\def \oint{{\cal O}_{int}}
\def \pos{{\mathbb{Z}}_{\geq 0}}
\def \xi{\Xi}

\def \dim{\mathrm{dim}}

\def \half{\frac{1}{2}}

\newcounter{dummy}
\theoremstyle{plain}
 
	 	\newtheorem{thm}{Theorem}[section]
	 	\newtheorem{lem}[thm]{Lemma}
	 	\newtheorem{cor}[thm]{Corollary}
		\newtheorem{prop}[thm]{Proposition}

\theoremstyle{definition}
 		\newtheorem{defn}[thm]{Definition}
		\newtheorem{rem}[thm]{Remark}
		\newtheorem{ntn}[thm]{Notation}
		\newtheorem{ex}[thm]{Example}

\begin{document}
\thispagestyle{empty}

\vskip 1.5 cm

\begin{center}
	\noindent{\Large \textbf{%
		Crystal Bases for the Quantum Superalgebra 
	}}\\

	\vspace*{0.5cm}

	\noindent{\Large \textbf{%
		{\boldmath$U_q(D(N,1))$}, {\boldmath$U_q(B(N,1))$}
	}}\\
	\renewcommand{\thefootnote}{\fnsymbol{footnote}}

	\vskip 2cm

	{\LARGE
		SUZUKI Kenei%
	}\\


	\noindent{ \bigskip }\\

	\it
		Graduate School of Mathematical Sciences, University of Tokyo \\
		Komaba 3-8-1, Meguro-ku, Tokyo 153-8914, Japan
	\noindent{\smallskip  }\\
		e-mail address: kenei@ms557pur.ms.u-tokyo.ac.jp%

	\bigskip
\end{center}
\begin{abstract}

 Let $V(\lam)$ be the irreducible lowest weight $U_q(D(N,1))$-module 
 with lowest weight $\lam$.
 Assume $\lam = n_0\ome0-\sum_{i=0}^{N}n_i\omg_i$, $n_i \in \pos$,
 where $\ome0$ is the fundamental weight corresponding to
 the unique odd coroot $h_0$.
 $V(\lam)$ is called typical if $n_0 \geq 0$.
 In this paper, we construct polarizable crystal bases of $V(\lam)$
 in the category $\oint$, which is a class of 
 integrable modules. 
 We also describe the decomposition of the tensor product of typical representations
 into irreducible ones, using the generalized Littlewood-Richardson rule for $\udn$.

 We also present analogous results for the quantum superalgebra $U_q(B(N,1))$.
\end{abstract}	

 {\it keywords:} Crystal bases; Quantum superalgebras; $U_q(D(N,1))$;
 $U_q(B(N,1))$; Tensor products

\vfill

\vfill
\setcounter{footnote}{0}
\renewcommand{\thefootnote}{\arabic{footnote}}


%
%

\section{Introduction}
\label{sec:intro}

        The theory of crystal base for quantized Lie algebras initiated by
	Kashiwara has brought and is still bringing a great deal of fruits
	in representation theory.
	It is hence natural to generalize it to the
	case of Lie superalgebras. 
	
        Finite dimensional
        simple Lie superalgebras $\gge$ were classified by Kac \cite{Kac1}. 
	He also determined the conditions under which a $\gge$-module
	is finite dimensional.
	The quantized enveloping algebras for Lie superalgebras are
	defined 
	by Yamane \cite{Yamane} for finite dimensional contragredient Lie superalgebras,
	and by Benkart, Kang and Melville \cite{BKD} for Borcherds
	superalgebras.
	Borcherds superalgebras include Kac-Moody superalgebras, which are
	analogous to Kac-Moody Lie algebras.

	The first work on the crystal bases for quantized Lie
	superalgebras was the one by Musson and Zou \cite{MZ}.
	They defined and constructed crystal bases of finite
	dimensional modules of 
	$U_q(B(0,N))=U_q({\mathfrak{osp}}(1,2N))$. 
	The crystal bases for $U_q(B(0,N))$ are essentially 
	the same as those for $U_q(B(N))$.
	Among the Lie superalgebras in Kac's list,
	$B(0,N)$ is distinguished because it is the unique one
	whose finite dimensional representations are completely reducible.
	$B(0,N)$ is also a Kac-Moody superalgebra, 
	while other finite dimensional simple Lie superalgebras are not.
	Jeong \cite{J} generalized these results from the point of view
	of Kac-Moody superalgebras.
	He defined crystal bases and showed their existence for the
	quantized Kac-Moody superalgebras when the modules are
	integrable. 


	Benkart, Kang, Kashiwara \cite{BKK} defined
        crystal bases for quantum contragredient superalgebras which are not 
	Kac-Moody Lie superalgebras
	(i.e. the ones containing $\otimes$ in their Dynkin diagrams).
	They introduced 
	a category $\oint$ of $\uge$-modules for
	contragredient Lie superalgebra $\gge$ (see Definition
        \ref{oint}),
	and defined the notion of crystal base for modules $M$ in $\oint$.
	The Kashiwara operators for odd roots behave quite differently
	from the case of Lie algebras 
	and the Kac-Moody Lie superalgebras 
	(see (\ref{Kashi0}) and Proposition \ref{tensor rule}).
	They also showed the existence of polarizable crystal bases
	of finite dimensional irreducible
	$U_q(\mathfrak{gl}(m,n))$-modules under some conditions on
	the highest weight and described them in terms of Young tableaux.
	Benkart and Kang \cite{BK} give a concise review for these
	results for quantized Lie superalgebras.

	A noteworthy feature of $U_q(\mathfrak{gl}(m,n))$ is that
	the vector representation belongs to $\oint$.
	However
	this fails for the other types of classical Lie superalgebras
	$B(m,n)$ $(m \geq 1)$,
	$C(n),D(m,n)$ and $D(2,1;\alpha)$.
	Therefore, when one attempts to generalize the results of
	\cite{BKK}, the first question is to find modules in $\oint$.
	
	The first generalization to these quantum superalgebras was for $U_q(D(2,1;\alpha))$
        by Zou \cite{Z1}.
	He found highest weight modules in $\oint$
	which are infinite dimensional,
        and constructed crystal bases of them. 
	In \cite{Z1}, $\alpha$ is assumed to be 
	an integer satisfying $\alpha \leq -2$.
	The reason he studied infinite dimensional modules is that any finite
        dimensional $U_q(D(2,1;\alpha))$-module $M$ does not satisfy the condition (iv)
	in the definition of $\oint$, that is 
	\begin{equation}\label{cond4}
	 \mbox{(iv)}\quad
	  \mbox{For any } \mu \in P,\quad M_\mu \neq 0  \mbox{ implies } 
					\langle h_0,\mu\rangle \ge 0,
	\end{equation}
	where $h_0$ is the unique simple odd coroot of $U_q(D(2,1;\alpha))$.

	\bigskip

	In this paper, we present two results on crystal bases for 
	$U_q(D(N,1))$. 
	Let $V(\lam)$ be the irreducible lowest
	weight module with lowest weight 
	\begin{equation}\label{form}
		\lambda = n_0{\omega}_0 - \sum_{i=1}^{N}n_i{\omega}_i, \quad n_i \in \pos \quad\mbox{ for }
		 0 \leq i \leq N,
	\end{equation}
	where $\omg_i$ are the fundamental weights.
	The first result is that 
	$V(\lam)$ admits a crystal base $B(\lam)$
	 (Theorem \ref{exist}).
	The weight (\ref{form}) is said to be typical if $n_0 \geq 1$.
	As the second result, we describe the decomposition of the tensor product of 
	the irreducible representations which have typical lowest weights (Theorem \ref{main}). 
	Our description heavily relies on the
	generalized Littlewood-Richardson rule for $\udn$ by Nakashima \cite{Nakashima}.
	We give similar results also for the algebra $U_q(B(N,1))$.

	Our work is a kind of generalization of Zou's work,
	since $U_q(D(2,1;\alpha)) = U_q(D(2,1))$ holds if $\alpha = 1$ and
	there is an isomorphism of algebras between 
	$U_q(D(2,1;\alpha))$ and $U_q(D(2,1;-1-\alpha))$.
	However, it should be noted that
	we have to adopt another approach because
	a $U_q(D(2,1;\alpha))$-module in $\oint$  
	does not necessarily belong to $\oint$
	when viewed as a $U_q(D(2,1;-1-\alpha))$-module.

	As in the case of $U_q(D(2,1;\alpha))$, (1.1) fails for any finite dimensional
	$U_q(D(N,1))$-module. 
	In addition,
	infinite dimensional highest weight modules do not belong to $\oint$.
	We first show that the irreducible lowest weight module
	$V(-\omg_N)$ with lowest weight $-\omg_N$ belongs to $\oint$,
	and give an explicit construction of a crystal base 
	$B(-\omg_N)$ (Proposition \ref{superspin}).
 	We find that $B(-\omg_N)$ is indexed by the crystal bases of the spin representations
	for $\udn$ (denoted by $\ovB_{sp}^{\pm}$ in this paper)
	and non-negative integers.
	Next we show that $B(-\omg_N) \otimes B(-\omg_N)$ contains
	$B(\ome0),B(-\omg_1),\ldots,B(-\omg_N)$
	using the decomposition of $\ovB_{sp}^+ \otimes \ovB_{sp}^{\pm}$.
	The existence of the crystal base $B(\lam)$ for general $\lam$
	(\ref{form}) follows
	by taking the tensor product of them.

	If a weight $\lam'$ is typical, we observe that
	$b'\otimes b \in B(\lam') \otimes B(\lam)$
	is the lowest weight vector for $U_q(D(N,1))$
	if and only if 
	$b'\otimes b \in B(\lam') \otimes B(\lam)$
	is the lowest weight vector for $\udn$.
	We use this fact to study the tensor product of typical representations.
	To get the lowest weight vectors of this tensor product,
	we decompose $B(\lam)$
	into copies of crystal bases of $\udn$ 
	labeled by certain non-negative integers
	(Proposition \ref{lamome0}).
	The typicality of $B(\lam)$ and the generalized
	Littlewood-Richardson rule enable us to find this decomposition.
	The integer labels are obtained through properties on
	Young tableaux for $\udn$ due to Koga \cite{Koga} 
	(see Proposition \ref{Koga} and Lemma \ref{0connection}).
	Applying the generalized Littlewood-Richardson rule again,
	we obtain all the lowest weight vectors of $B(\lam') \otimes B(\lam)$.
	This is Theorem \ref{main}.
	
	\bigskip

	This paper is organized as follows.
	In Section 2, we give the definition of the quantized Lie superalgebra 
	$U_q(D(N,1))$ following Yamane \cite{Yamane} and set up the notation.
	In Section 3, we review basic facts about the category $\oint$ and the crystal bases
	in the context of $U_q(D(N,1))$.
%
        In Section 4, we recall some properties of the crystal bases for
        $\udn$ which we use in this paper.	
	In Section 5, we construct a crystal base of $V(-\omg_N)$.
        In Section 6 and Section 7, we state the existence of crystal
	bases and give a decomposition of $B(\lam)$ mentioned above for
	typical $\lam$.
	In Section 8, we decompose the tensor product of 
	typical representations into irreducible ones in Theorem \ref{main}.
	We also present examples for Theorem \ref{main}
	in the case of $U_q(D(2,1))$ and $U_q(D(4,1))$.
	The results for $U_q(B(N,1))$ are summarized in Section 9.

\section{Definition of {\boldmath$U_q(D(N,1))$}}

	In this section we fix our notation concerning the quantum universal enveloping
        algebra for the Lie superalgebra $\gge = D(N,1)$.

	Let $P=\oplus_{i=0}^{N}{\mathbb{Z}}\omg_i$ be a free
        $\Z$-module with basis $\{\omg_i\}_{i=0}^{N}$.
	Let $\{h_i\}_{i=0}^{N}$ be the dual basis of $P^* =
        \mathrm{Hom}_{\mathbb Z} (P,\mathbb Z)$ (simple
        coroots) relative
        to the pairing $\left<\cdot,\cdot\right>$.
	Define simple roots $\{\alpha_i\}_{i=0}^{N} \subset P$ so
        that $a_{ij} = \left<h_i,\alpha_i\right>$ is given by
        following Cartan matrix $A=(a_{ij})_{0 \leq i,j \leq N}$;

	\begin{eqnarray}
	  \left(\begin{array}{rrrrrrrr}\nonumber

	    0       & 1      & 0      & \cdots & \\

	    -1      & 2      & -1     &\\
	    0       & -1     & 2      & -1     &\\
	    \vdots  &        & \ddots & \ddots & \ddots&       & \vdots & \vdots\\
                    &        &        &  -1    & 2     & -1    & 0      &0\\
                    &        &        &        & -1    & 2     & -1     &-1\\
	            &        &        & \cdots & 0     & -1    & 2      &0\\
                    &        &        & \cdots & 0     & -1    & 0      &2\\
	                       \end{array}\right).
	\end{eqnarray}

	The associated Dynkin diagram is 

	\begin{center}
\unitlength 0.1in
\begin{picture}( 44.5000, 10.0000)(  5.0000,-11.0000)
%
\special{pn 8}%
\special{ar 600 600 100 100  0.0000000 6.2831853}%
%
\special{pn 8}%
\special{ar 1200 600 100 100  0.0000000 6.2831853}%
%
\special{pn 8}%
\special{ar 1800 600 100 100  0.0000000 6.2831853}%
\put(5.5000,-10.0000){\makebox(0,0)[lb]{$0$}}%
\put(11.5000,-10.0000){\makebox(0,0)[lb]{$1$}}%
\put(17.5000,-10.0000){\makebox(0,0)[lb]{$2$}}%
\put(49.5000,-2.7000){\makebox(0,0)[lb]{$N-1$}}%
\put(49.5000,-10.7000){\makebox(0,0)[lb]{$N$}}%
\put(39.6000,-9.4000){\makebox(0,0)[lb]{$N-2$}}%
%
\special{pn 8}%
\special{sh 1}%
\special{ar 2400 600 10 10 0  6.28318530717959E+0000}%
\special{sh 1}%
\special{ar 2600 600 10 10 0  6.28318530717959E+0000}%
\special{sh 1}%
\special{ar 2800 600 10 10 0  6.28318530717959E+0000}%
\special{sh 1}%
\special{ar 3000 600 10 10 0  6.28318530717959E+0000}%
\special{sh 1}%
\special{ar 3200 600 10 10 0  6.28318530717959E+0000}%
\special{sh 1}%
\special{ar 3400 600 10 10 0  6.28318530717959E+0000}%
\special{sh 1}%
\special{ar 3600 600 10 10 0  6.28318530717959E+0000}%
\special{sh 1}%
\special{ar 3600 600 10 10 0  6.28318530717959E+0000}%
%
\special{pn 8}%
\special{pa 530 530}%
\special{pa 670 670}%
\special{fp}%
\special{pa 670 530}%
\special{pa 530 670}%
\special{fp}%
%
\special{pn 8}%
\special{pa 700 600}%
\special{pa 1100 600}%
\special{fp}%
\special{pa 1300 600}%
\special{pa 1700 600}%
\special{fp}%
\special{pa 1900 600}%
\special{pa 2300 600}%
\special{fp}%
\special{pa 4100 600}%
\special{pa 3700 600}%
\special{fp}%
\special{pa 4200 600}%
\special{pa 4800 200}%
\special{fp}%
\special{pa 4800 1000}%
\special{pa 4200 600}%
\special{fp}%
%
\special{pn 8}%
\special{sh 0}%
\special{ar 4800 200 100 100  0.0000000 6.2831853}%
%
\special{pn 8}%
\special{sh 0}%
\special{ar 4200 600 100 100  0.0000000 6.2831853}%
%
\special{pn 8}%
\special{sh 0}%
\special{ar 4800 1000 100 100  0.0000000 6.2831853}%
\end{picture}%
 \mbox{ . }
	\end{center}

	We put
	\begin{eqnarray}\label{li}
	l_0 = 1,\quad l_1 = \cdots = l_N     = -1, 
	\end{eqnarray}
	and introduce a symmetric bilinear form $(\cdot,\cdot)$ on
	$\mathfrak{h}^*=P \otimes_{\mathbb{Z}} {\mathbb{C}}$ by
	\begin{equation}\nonumber
	l_i\left<h_i,\lam\right>=(\alpha_i,\lam) \quad\mbox{for any}\quad \lam \in P,\quad
	0 \leq i \leq N.
	\end{equation}

	\begin{defn}[Yamane \cite{Yamane} Theorem 10.5.1 ]\label{yama}
		Let $U'_q(D(N,1))$ be the associative algebra over
	        ${\mathbb{Q}}(q)$ with $1$ 
		generated by $e_i$, $f_i$, $q^h$ ($0 \leq i \leq N$, $h \in P^*$) with the following relations.
		\begin{equation}
		  q^0 = 1,\quad q^hq^{h'} = q^{h+h'} \quad\mbox{ for } h,h' \in P^*,
		\end{equation}
		\begin{equation}
		  q^he_iq^{-h} = q^{\langle h,{\alpha}_i \rangle}e_i,\quad
		  q^hf_iq^{-h} = q^{\langle -h,{\alpha}_i \rangle}f_i
		  \quad\mbox{ for } 0 \leq i \leq N \mbox{,}h \in P^*,
		\end{equation}
		\begin{equation}\label{efcom}
		  e_if_j - (-1)^{p(i)p(j)}f_je_i = 
		  \delta_{ij}\frac{t_i - {t_i}^{-1}}{q_i - {q_i}^{-1}}\quad\mbox{ for } 0 \leq i,j \leq N,
		\end{equation}
		\begin{equation}\label{commutation}
		 \begin{split}
		  \sum_{\nu = 0}^{1 + |a_{ij}|} {(-1)^\nu} \left[ \begin{array}{c}
		    1 + |a_{ij}|\\
		    \nu
		    \end{array} \right]_{q_i}
		  e_i^{1+|a_{ij}|-{\nu}}e_je_i^{\nu} =
		  \sum_{\nu = 0}^{1 + |a_{ij}|} {(-1)^\nu} \left[ \begin{array}{c}
		    1 + |a_{ij}|\\
		    \nu
		    \end{array} \right]_{q_i}
		  f_i^{1+|a_{ij}|-{\nu}}f_jf_i^{\nu} = 0, \\
		  \hfill \mbox{for }  1 \leq i \leq N , 0 \leq j \leq N ,i \neq j,
		 \end{split}
		\end{equation}

		\begin{equation}\label{square}
		  e_0^2=0 \mbox{,}\quad f_0^2=0,
		\end{equation}
	        where
	        \begin{eqnarray}\nonumber
		 t_i = q^{l_ih_i}, 
		\end{eqnarray}
	        \begin{equation}\nonumber
		 p(0) = 1,\quad p(i) = 0 \mbox{ for }  1 \leq i \leq N,
		\end{equation}
	        \begin{equation}\nonumber
		  \left[ \begin{array}{c}
		    m\\
		    n
		    \end{array} \right]_t=
		  \prod_{i=0}^{n-1}\frac{t^{m-i}-t^{-(m-i)}}{t^{i+1}-t^{-(i+1)}}.
		\end{equation}
	 The algebra $U_q'(D(N,1))$ is ${\mathbb Z}_2$-graded. Define the
	 parity operator $\sigma$ by setting
	 $\sigma(e_i)=(-1)^{p(i)}e_i$,
	 $\sigma(f_i)=(-1)^{p(i)}f_i$,
	 $\sigma(q^h)= q^h$ ($h \in P$).
%
	 Then the quantized Lie superalgebra $U_q(D(N,1))$ is
	 defined to be $U_q(D(N,1)) = U'_q(D(N,1)) \oplus U'_q(D(N,1)) \sigma$ with
	 the algebra structure given by $\sigma^2 = 1 $ and $\sigma u \sigma = \sigma(u)$ for 
	 $u \in U'_q(D(N,1))$.
	\end{defn}

	The Hopf algebra structure on $U_q(D(N,1))$ is given as follows.\\
	The comultiplication $\Delta$, the antipode $S$, and the counit $\varepsilon$ are defined by
	\begin{equation}\nonumber
	 \begin{array}{lll}
	  \Delta(\sigma) = \sigma \otimes \sigma,&     
	   S(\sigma) = \sigma,&
	    \varepsilon(\sigma)=1,\\
	  \Delta(q^h)    = q^h    \otimes q^h,&
	   S(q^h)    = q^{-h},&
	    \varepsilon(q^h)=1,\\
	  \Delta(e_i)    = e_i    \otimes t_i^{-1} + \sigma^{p(i)} \otimes e_i,&
	   S(e_i)    = -\sigma^{p(i)}e_it_i,&
	    \varepsilon(e_i)=0,\\
	  \Delta(f_i)    = f_i    \otimes 1        + \sigma^{p(i)}t_i \otimes f_i,&
	   S(f_i)    = -\sigma^{p(i)}t_i^{-1}f_i,&
	    \varepsilon(f_i)=0.\\
	 \end{array}
	\end{equation}

	We define an orthogonal basis 
	$\{$ $\delta$, $\varepsilon_1$,$\ldots$,$\varepsilon_N$ $\}$ 
	 of $\mathfrak{h}^*$ by

	\begin{equation}\nonumber
	 \alpha_0=\delta-\varepsilon_1,\quad \alpha_1=\varepsilon_1-\varepsilon_2,\ldots,
	 \alpha_{N-1}=\varepsilon_{N-1}-\varepsilon_N, \quad \alpha_N=\varepsilon_{N-1}+\varepsilon_N.
	\end{equation}
	We have
	\begin{equation}\nonumber
	 (\delta,\delta)=1,\quad (\varepsilon_i,\varepsilon_i)=-1 \quad\mbox{for } 1 \leq i \leq N,
	\end{equation}
	\begin{align*}
	 \omg_0&=\delta,\\
	 \omg_i&=-\delta + \epsilon_1 + \cdots + \epsilon_i \quad \mbox{for }1 \leq i \leq N-2,\\
	 \omg_{N-1}&=\half(-\delta + \epsilon_1 + \cdots + \epsilon_{N-1}- \epsilon_N),\\
	 \omg_{N}  &=\half(-\delta + \epsilon_1 + \cdots + \epsilon_{N-1}+ \epsilon_N). 
	\end{align*}
	The set of even and odd roots are 
	$\Delta_i = \Delta_i^+ \cup (-\Delta_i^+)$ ($i=0,1$), where
	\begin{equation}\nonumber
	 \Delta_0^+=\{\varepsilon_i \pm \varepsilon_j,2\delta\}_{1 \leq i < j \leq N}, \quad 
	 \Delta_1^+=\{\pm\varepsilon_i + \delta\}_{1 \leq i \leq N}.
	\end{equation}
	We put $\overline{\Delta_1}=\Delta_1$ and 
	\begin{equation}\nonumber
	 \rho_i = \sum_{\beta \in \Delta_i^+}\beta \quad \mbox{for } i=0,1, \qquad \rho = \rho_0 - \rho_1.
	\end{equation}

	We denote by $V(\lam)$ the irreducible lowest weight module with
	lowest weight
	$\lam$.
	\begin{defn}[Kac \cite{Kac2} Theorem 1]
	 For $\lam \in P$, $V(\lambda)$ is
	 a typical representation if
	 \begin{eqnarray}\nonumber
	  (\lambda - \rho , \beta) \neq 0 \mbox{ for any } \beta \in \overline{\Delta_1}.
	 \end{eqnarray}
	 In this case, $\lambda$ is called a typical weight. A weight
	 which is not typical is called an atypical weight.
	\end{defn}
	 In this paper we will study $V(\lam)$ where
	 $\lambda = n_0\omega_0 - \sum_{i=1}^{N}n_i\omega_i$ with
	 $n_i \in {\mathbb{Z}}_{\geq0}$ for $0 \leq i \leq N$. It is typical if and
	 only if $n_0 \geq 1$.


\section{Category {\boldmath${\cal O}_{int}$} and Crystal Base}

	\begin{defn}[\cite{BKK} Definition 2.2]\label{oint}
		$\oint$ is the category whose objects are $U_q(D(N,1))$-modules $M$ and 
		whose morphisms are $U_q(D(N,1))$-linear homomorphisms satisfying the following conditions:

			\begin{description}
				\item[(i)]
					$M$ has a weight decomposition 
					   $M=\bigoplus_{\lambda \in P}M_\lambda$, where\\
					$M_\lambda = \{u \in M ; q^{h}u = q^{\left< h,\lambda \right>}u 
					\mbox{ for any } h \in P^* \}$,

				\item[(ii)]
					$\dim M_\lambda < \infty$ for any $\lambda \in P$, 
				
				\item[(iii)]
					For $1 \leq i \leq N$, $M$ is  locally $U_q(D(N,1))_i$-finite,
					that is, $\dim U_q(D(N,1))_i u < \infty$ for  any $u \in M$, 
				\item[(iv)]
					For any $\lam \in P$, $M_\lam \neq 0 $ implies 
					$\langle h_0,\lam\rangle \ge 0$,

				\item[(v)]
					$e_0M_\lam = f_0M_\lam = 0$ for any $\lam \in P $ such that 
					$\langle h_0,\lam \rangle = 0$.

			\end{description}

	\end{defn}

	We define the Kashiwara operators $\tei$, $\tfi$ on $M \in \oint$. 
	For any
	$u \in M_{\lam}$, and $i={1,\ldots,N}$, let
	\begin{eqnarray}\nonumber
	 u = \sum_{k \geq 0,-\langle h_i, \lambda \rangle}f_i^{(k)}u_k
	\end{eqnarray}
	be the unique expression with $e_iu_k=0$ for each $k$.
	We define
	\begin{eqnarray}\nonumber
	 \begin{split}
	 \widetilde{e_i}u &= \sum_kf_i^{(k-1)}u_k,\\
	 \widetilde{f_i}u &= \sum_kf_i^{(k+1)}u_k.
	 \end{split}
	\end{eqnarray}
	For $i=0$, define
	\begin{align}\label{Kashi0}
	 \widetilde{e_0}u &= q_0^{-1}t_0e_0u,\\
	 \widetilde{f_0}u &= f_0u.
	\end{align}
	
	We set
	\begin{equation}\nonumber
	 A=\left\{ \frac{f}{g};f,g\in {\mathbb{Q}}[q] , g(0) \neq 0 \right\}.
	\end{equation}

	\begin{defn}[\cite{BKK} Definition 2.3, 2.4]
	 Let $M$ be a $U_q(D(N,1))$-module in the category
	 ${\cal O}_{int}$. A crystal base of $M$ is a pair $(L,B)$ such that
	 \begin{description}
	  \item[(L1)]
		     $L$ is a free $A$-submodule satisfying $M = L \otimes_A {\mathbb{Q}}(q)$,
	  \item[(L2)]
		     $\sigma L = L$ and $L$ has a weight decomposition
		     $L = \oplus_{\lambda \in P}L_{\lambda}$ with
		     $L_{\lambda} = L \cap M_{\lambda}$,
          \item[(L3)]
		     $\widetilde{e_i}L \subset L$ and $\widetilde{f_i}L \subset L$ for $0 \leq i \leq N$,
          \item[(B1)]
		     $B$ is a subset of $L/qL$ such that $\sigma b = \pm b$ for any $b \in B$, and 
		     $B = \sqcup_{\lambda \in P}B_{\lambda}$ with 
		     $B_{\lambda} = B \cap (L_{\lambda}/{qL_{\lambda}})$,
	  \item[(B2)]
		     $B$ is a pseudo-base of $L/{qL}$, that is , $B = B' \cup (-B')$ for a $\mathbb{Q}$
		     basis $B'$ of $L/qL$,
          \item[(B3)]
		     $\widetilde{e_i}B \subset B \sqcup \{0\}$ and
		     $\widetilde{f_i}B \subset B \sqcup \{0\}$ for $0 \leq i \leq N$,
          \item[(B4)]
		     For any $b,b' \in B$, and $0 \leq i \leq N$, 
		     $b = \widetilde{f_i}b' \mbox{ if and only if } \widetilde{e_i}b = b'$.

	 \end{description}
	\end{defn}
	We often denote $(L,B)$ by $B$. 

	Let $\eta$ be the anti-automorphism of $U_q(D(N,1))$ defined by
	\begin{equation}\nonumber
	 \eta(\sigma)=\sigma,\quad \eta(q^h)=q^h, \quad \eta(e_i)=q_if_it_i^{-1}, \quad \eta(f_i)=q_i^{-1}t_ie_i.
	\end{equation}
	A symmetric bilinear form $(\cdot,\cdot)$ on a $U_q(D(N,1))$-module $M$
	is called polarization if $(au,v)=(u,\eta(a)v)$ holds for any 
	$u,v \in M$ and $a \in U_q(D(N,1))$.
	\begin{defn}
	 A crystal base $(L,B)$ for a $U_q(D(N,1))$-module $M$ is polarizable
	 if there exists a polarization $(\cdot,\cdot)$ of $M$ such that
	 $(L,L) \subset A$, and with respect to the induced
	 ${\mathbb{Q}}$-valued symmetric bilinear form $(\cdot,\cdot)_0$ on
	 $L/qL$,
	 \begin{equation}\nonumber
	  (b,b')_0 = \begin{cases}\pm1 \quad &\mbox{if }b'=\pm b, \\
		                  0          &\mbox{otherwise}
		     \end{cases}
	 \end{equation}
	 holds for any $b$, $b' \in B$.
	\end{defn}

	\begin{rem}\label{cond4}
	Theorem 8 of \cite{Kac1} implies that condition {\bf{(iv)}} of
	 Definition \ref{oint} fails for
	 all non-trivial irreducible 
	finite dimensional $U_q(D(N,1))$-modules.
	For this reason, we treat infinite dimensional lowest weight modules.
	\end{rem}

	For $b \in B$, $0 \leq i \leq N$, we define
	\begin{equation}\nonumber
	 \varepsilon_i(b)=\max\{n \in \pos ; \tei^nb \neq 0\},\quad
	 \varphi_i(b)    =\max\{n \in \pos ; \tfi^nb \neq 0\}.
	\end{equation}
	We write $wt(b) = \lam$ for $b \in B_{\lam}$.

	\begin{prop}[\cite{BKK} Proposition 2.8]\label{tensor rule}
	  Let $(L_i,B_i)$ be polarizablue crystal bases of
	 $U_q(D(N,1))$-modules $M_i \in \oint$, $i=1,2$. Then
	  $(L_1 \otimes_A L_2,B_1 \otimes B_2)$ is a polarizable crystal base of $M_1 \otimes M_2$, and
	  the actions of $\widetilde{e_i}$ and $\widetilde{f_i}$ are given as follows.\\
	 \begin{description}
	  \item[(1)]$i=1,2,\dots,N$
	  \begin{eqnarray}\label{etildei}
	    \widetilde{e_i}(b_1 \otimes b_2) = \begin{cases}
                                                  b_1 \otimes \widetilde{e_i}(b_2) &\text{ if } 
						   \varepsilon_i(b_1) \leq \varphi_i(b_2),\\
	                                          \widetilde{e_i}(b_1) \otimes b_2 &\mbox{ if } 
						    \varepsilon_i(b_1) > \varphi_i(b_2),
						  \end{cases}
	  \end{eqnarray}
	  \begin{eqnarray}\label{ftildei}
	    \widetilde{f_i}(b_1 \otimes b_2) = \begin{cases}
                                                  b_1 \otimes \widetilde{f_i}(b_2) &\text{ if } 
						    \varepsilon_i(b_1) < \varphi_i(b_2),\\
	                                          \widetilde{f_i}(b_1) \otimes b_2 &\mbox{ if }
						    \varepsilon_i(b_1) \geq \varphi_i(b_2),
						  \end{cases}
	  \end{eqnarray}

	  \item[(2)]$i=0$
	  	  \begin{eqnarray}\label{etilde0}
	                \widetilde{e_0}(b_1 \otimes b_2) = \begin{cases}
                                                  \sigma b_1 \otimes \widetilde{e_0}(b_2) &\text{ if } 
	                                           \langle h_0,wt(b_1)\rangle = 0, \\
	                                          \widetilde{e_0}(b_1) \otimes b_2 &\text{ if } 
                                                   \langle h_0,wt(b_1)\rangle > 0,
						  \end{cases}
	  \end{eqnarray}
    	  \begin{eqnarray}\label{ftilde0}
	                \widetilde{f_0}(b_1 \otimes b_2) = \begin{cases}
                                                  \sigma b_1 \otimes \widetilde{f_0}(b_2) &\text{ if } 
	                                           \langle h_0,wt(b_1)\rangle = 0, \\
	                                          \widetilde{f_0}(b_1) \otimes b_2 &\mbox{ if } 
                                                   \langle h_0,wt(b_1)\rangle > 0.
						  \end{cases}
	  \end{eqnarray}
         \end{description}
	\end{prop}

	Note that in (\ref{etildei})-(\ref{ftildei}) the inequality
	signs are
	opposite to those for ordinary Lie algebra case.
	This is due to the negative sign of $l_i$.

	Let $(L,B)$ be a crystal base for a $U_q(D(N,1))$-module in
	$\oint$. We define
	\begin{equation}\nonumber
	 \begin{split}
	 LW(B) = \{b \in B \mbox{ }; \tfi(b) = 0 \mbox{ for } 0 \leq i \leq N\},\\
	 \overline{LW}(B) = \{b \in B \mbox{ }; \tfi(b) = 0 \mbox{ for } 1 \leq i \leq N\}.
	 \end{split}
	\end{equation}

        \begin{lem}\label{lwv}
	  Let $(L_i,B_i)$ be as in Proposition \ref{tensor rule}. Assume
	  $LW(B_i) \neq \emptyset $ for $i=1,2$. Then
	  the element of $LW(B_1 \otimes B_2)$ is of the form $u_1 \otimes v$
	  with $u_1 \in LW(B_1)$.

	  \begin{proof}
	    Assume that $u \otimes v \in LW(B_1 \otimes B_2)$ with $u \not\in LW(B_1)$. We have two cases;\\
	    Case 1: $\widetilde{f_i}u \neq 0 $ for some $1 \leq i \leq N$.\\
	    Since $0 = \widetilde{f_i}(u \otimes v) = u \otimes \widetilde{f_i}(v)$, (\ref{ftildei}) implies
	    $\varepsilon_i(u) < \varphi_i(v)$. Hence $\widetilde{f_i}(v) \neq 0$. This is a contradiction.\\
	    Case 2: $\widetilde{f_0}u \neq 0$, $\widetilde{f_i}u = 0$ for $1 \leq i \leq N$.\\
	    Since $0=\widetilde{f_0}(u \otimes v) = u \otimes \widetilde{f_0}(v)$, (\ref{ftilde0}) implies
	    $\langle h_0 , wt(u) \rangle = 0$. This contradicts
	    Definition \ref{oint}(v). 
	  \end{proof}
	\end{lem}

	In \cite{Z1}, the algebra $U_q(D(2,1;\alpha))$ with 
	$\alpha \leq -2$ is considered. Lemma
	\ref{lwv} fails in this case since the $l_i$ have
	both positive and negative signs. 

\section{Results on crystal bases for {\boldmath$\udn$}-modules}

	 The even part of $U_q'(D(N,1))$ is the eigenspace of $\sigma$
	 with eigenvalue $+1$,
	 denoted by $U_q'(D(N,1)_0)$.
	 In our case it is given by
	 $U_q(D(N)) \otimes U_q(C(1))$,
	 where
	 $U_q(D(N))$ is the subalgebra with generators 
	 $e_i, f_i,q^{h_i}$ $(1 \leq i \leq N)$,
	 and $U_q(C(1)) \simeq U_q(\mathfrak{sl}_2)$ is the one generated by
	 $E,F,q^H$, where $H=2(h_0-h_1-\cdots-h_{N-2})-h_{N-1}-h_N$,
	 $E$ and $F$ are the elements corresponding to the root
	 $2(\alpha_0 + \cdots + \alpha_{N-2}) + \alpha_{N-1} + \alpha_N$.
%
	In this section we recall known properties of
	crystal bases for $\udn$.

	We denote the irreducible lowest weight module of
	$U_q(D(N))$ with
	lowest weight $\Lambda$ and its crystal base
	by $\ovV(\Lam)$ and $\ovB(\Lambda)$ respectively,
	and refer to the crystal base as
	a $\udn$-crystal for short. 
	Let $\{\Lam_i\}_{i=1}^{N}$ be the fundamental weights of $U_q(D(N))$.
	We denote 
	$\ovB(-\Lam_{N-1})$ by
	$\ovB_{sp}^-$, and $\ovB(-\Lam_N)$ by $\ovB_{sp}^+$.
	
	The crystal bases of spin representations are realized as
	\begin{equation}\nonumber
	 \ovB_{sp}^+ \cong \{b=(i_1,\ldots,i_N) ; i_1,\ldots,i_N=\pm,
	  \mbox{ and } - \mbox{ appears an even number of times}\},
	\end{equation}
	\begin{equation}\nonumber
	 \ovB_{sp}^- \cong \{b=(i_1,\ldots,i_N) ; i_1,\ldots,i_N=\pm,
	  \mbox{ and } -  \mbox{ appears an odd number of times}\},
	\end{equation}
	with the lowest weight vectors $(+,\ldots,+)$ and $(+,\ldots,+,-)$ respectively.
	The actions of $\tei$ and $\tfi$ read
	\begin{eqnarray}\nonumber
	 \hspace*{10mm}
	    \widetilde{f_l}(i_1,i_2,\ldots,i_N) &=&
		    \begin{cases}
		     (i_1,\ldots,\stackrel{l}{+},\stackrel{l+1}{-},\ldots,i_N) &\mbox{ if } i_l=-,i_{l+1}=+,\\\nonumber 
		     0                       &\mbox{ otherwise, }
		    \end{cases}\\\nonumber 
	    \widetilde{e_l}(i_1,i_2,\ldots,i_N) &=&
		    \begin{cases}
		     (i_1,\ldots,\stackrel{l}{-},\stackrel{l+1}{+},\ldots,i_N) &\mbox{ if } i_l=+,i_{l+1}=-,\\\nonumber 
		     0                       &\mbox{ otherwise, }
		    \end{cases}\\\nonumber 
	\end{eqnarray}
	   \mbox{for $1 \leq l \leq N-1$ and,}
	\begin{eqnarray}\nonumber
	    \widetilde{f_N}(i_1,i_2,\ldots,i_N) &=&
		    \begin{cases}
		     (i_1,\ldots,\stackrel{N-1}{+},\stackrel{N}{+}) &\mbox{ if } i_{N-1}=-,i_N=-,\\\nonumber 
		     0                      &\mbox{ otherwise, }
		    \end{cases}\\\nonumber 
	    \widetilde{e_N}(i_1,i_2,\ldots,i_N) &=&
		    \begin{cases}
		     (i_1,\ldots,\stackrel{N-1}{-},\stackrel{N}{-}) &\mbox{ if } i_{N-1}=+,i_N=+,\\\nonumber 
		     0                      &\mbox{ otherwise. }
		    \end{cases}
	 \end{eqnarray}

	Set
	\begin{equation}\nonumber
	 \xi_0=0, \quad \xi_i = \Lam_i \quad(1 \leq i \leq N-2), \quad \xi_{N-1}=\Lam_{N-1}+\Lam_N,\quad
	 {\xi_N}'=2\Lam_{N-1}, \quad \xi_N=2\Lam_N.
	\end{equation}
        \begin{prop}[Nakashima \cite{Nakashima}]\label{spinspin}
        We have the decomposition of crystals
	 \begin{eqnarray}\label{clpp}
	  \ovB_{sp}^+ \otimes \ovB_{sp}^+ = \bigoplus_{\substack{0 \leqslant k \leqslant N\\
	                                                         k \equiv N \mod 2\\}}
	   \ovB(-\xi_k),
	 \end{eqnarray}
	 \begin{eqnarray}\label{clpm}
	  \ovB_{sp}^+ \otimes \ovB_{sp}^- = \bigoplus_{\substack{0 \leqslant k \leqslant N-1\\
	                                                         k \equiv N-1 \mod 2\\}}
	   \ovB(-\xi_k).
	 \end{eqnarray}

	 For $0 \leq i \leq N$, the lowest weight vector
	 corresponding to the connected component $\ovB(-\xi_i)$ is 
	 \begin{eqnarray}\label{clhwv}
	  (+,\ldots,+) \otimes (\overbrace{+,\ldots,+}^{i},-,\ldots,-).	   
	 \end{eqnarray}
	 \end{prop}

	 Each connected component $\ovB(-\xi_k)$ can be given an explicit characterization.
	 For that purpose it is convenient to use an alternative
	 description of $\ovB_{sp}^{\pm}$ in terms of semi-standard
	 Young tableaux \cite{KN}.
	 Consider the set of letters $S = \{1,\ldots,N,\overline{N},\ldots,\overline{1}\}$. 
	 We introduce an ordering
	 $\prec$ on $S$ by
	 \begin{equation}\nonumber
	  1 \prec 2 \prec \cdots \prec N-1 \prec \begin{matrix} N \\ \overline{N} \end{matrix} \prec
	  \overline{N-1} \prec \cdots \prec \overline{2} \prec \overline{1}.
	 \end{equation} 
	 Then, there is an isomorphism of crystals 
	 \begin{equation}\label{des}
	  \ovB_{sp}^+ \sqcup \ovB_{sp}^- \cong \left\{ \begin{minipage}[c]{7mm}
\unitlength 0.1in
\begin{picture}(  2.0000,  8.0000)(  4.0000,-12.0000)
%
\special{pn 8}%
\special{pa 400 400}%
\special{pa 600 400}%
\special{fp}%
\special{pa 600 400}%
\special{pa 600 1200}%
\special{fp}%
\special{pa 600 1200}%
\special{pa 400 1200}%
\special{fp}%
\special{pa 400 1200}%
\special{pa 400 400}%
\special{fp}%
\special{pa 400 400}%
\special{pa 400 800}%
\special{fp}%
\put(4.3000,-5.8000){\makebox(0,0)[lb]{$a_1$}}%
\put(4.9000,-9.0000){\makebox(0,0)[lb]{$\vdots$}}%
\put(4.2000,-11.8000){\makebox(0,0)[lb]{$a_N$}}%
%
\special{pn 8}%
\special{pa 400 600}%
\special{pa 600 600}%
\special{fp}%
\special{pa 600 1000}%
\special{pa 400 1000}%
\special{fp}%
\end{picture}%

				  \end{minipage};
				  \begin{split}
				  &\mbox{ (1) } a_1,\ldots,a_N \in S\\
	                          &\mbox{ (2) } a_1 \prec \cdots \prec a_N\\
				  &\mbox{ (3) } a \mbox{ and }\overline{a} \mbox{ do not appear simultaneously}
				  \end{split}
			  \right\}.
	 \end{equation}
	 In this description, $a$ corresponds to the $a$-th $+$ and
	 $\overline{a}$ corresponds to the $a$-th $-$ in the former
	 description.
	 
	 We introduce notations of Young tableaux for convenience.

	 \begin{ntn}
	  \begin{description}
	   \item{(1)}
		       A Young tableau  \begin{minipage}[c]{7mm}
\unitlength 0.1in
\begin{picture}(  2.0000,  8.0000)(  4.0000,-12.0000)
%
\special{pn 8}%
\special{pa 400 400}%
\special{pa 600 400}%
\special{fp}%
\special{pa 600 400}%
\special{pa 600 1200}%
\special{fp}%
\special{pa 600 1200}%
\special{pa 400 1200}%
\special{fp}%
\special{pa 400 1200}%
\special{pa 400 400}%
\special{fp}%
\special{pa 400 400}%
\special{pa 400 800}%
\special{fp}%
\put(4.3000,-5.8000){\makebox(0,0)[lb]{$a_1$}}%
\put(4.9000,-9.0000){\makebox(0,0)[lb]{$\vdots$}}%
\put(4.2000,-11.8000){\makebox(0,0)[lb]{$a_N$}}%
%
\special{pn 8}%
\special{pa 400 600}%
\special{pa 600 600}%
\special{fp}%
\special{pa 600 1000}%
\special{pa 400 1000}%
\special{fp}%
\end{picture}%

				        \end{minipage}
		      is denoted by $t(a_1,\ldots,a_N)$.

	   \item{(2)}
		       A skew Young tableau \begin{minipage}[c]{28mm}
\unitlength 0.1in
\begin{picture}( 14.7000, 14.0000)( 2.5000,-18.0000)
\put(10.0000,-5.8000){\makebox(0,0)[lb]{$a_1$}}%
\put(8.1200,-9.8000){\makebox(0,0)[lb]{$a_{N-k+1}$}}%
\put(10.0000,-13.8000){\makebox(0,0)[lb]{$a_N$}}%
\put(6.0000,-9.8000){\makebox(0,0)[rb]{$b_1$}}%
\put(6.0000,-13.7000){\makebox(0,0)[rb]{$b_k$}}%
\put(6.2000,-17.7000){\makebox(0,0)[rb]{$b_N$}}%
%
\special{pn 8}%
\special{pa 800 400}%
\special{pa 1300 400}%
\special{fp}%
\special{pa 1300 400}%
\special{pa 1300 1400}%
\special{fp}%
\special{pa 1300 1400}%
\special{pa 300 1400}%
\special{fp}%
\special{pa 300 1400}%
\special{pa 300 800}%
\special{fp}%
\special{pa 300 800}%
\special{pa 1300 800}%
\special{fp}%
\special{pa 1300 1000}%
\special{pa 300 1000}%
\special{fp}%
\special{pa 300 1200}%
\special{pa 1300 1200}%
\special{fp}%
\special{pa 800 400}%
\special{pa 800 1800}%
\special{fp}%
\special{pa 800 1800}%
\special{pa 300 1800}%
\special{fp}%
\special{pa 300 1400}%
\special{pa 300 1800}%
\special{fp}%
\special{pa 300 1400}%
\special{pa 300 800}%
\special{fp}%
\special{pa 300 800}%
\special{pa 300 1400}%
\special{fp}%
\special{pa 300 1600}%
\special{pa 800 1600}%
\special{fp}%
\special{pa 800 600}%
\special{pa 1300 600}%
\special{fp}%
\put(10.0000,-7.7000){\makebox(0,0)[lb]{$\vdots$}}%
\put(10.0000,-11.7000){\makebox(0,0)[lb]{$\vdots$}}%
\put(5.5000,-11.7000){\makebox(0,0)[rb]{$\vdots$}}%
\put(5.5000,-15.7000){\makebox(0,0)[rb]{$\vdots$}}%
\end{picture}%

				            \end{minipage}
		      is denoted by\\ 
		      $t(a_1,\ldots,a_{N-k};a_{N-k+1},\ldots,a_N|b_1,\ldots,b_k;b_{k+1},\ldots,b_N)$.
		      
		       Note that in this skew Young tableau, the number of rows which has two boxes
		      is $k$.
	  \end{description}
	 \end{ntn}

	 \begin{defn}
	  A skew Young tableau
	  $t(a_1,\ldots,a_{N-k};a_{N-k+1},\ldots,a_N|b_1,\ldots,b_k;b_{k+1},\ldots,b_N)$
	  is semi-standard if
	  \begin{align*}
	   &a_1,\ldots,a_N \mbox{ satisfy (1),(2) and (3) in (\ref{des})},\\
	   &b_1,\ldots,b_N \mbox{ satisfy (1),(2) and (3) in (\ref{des})},\\
	   &b_r \preceq a_{N-k+r} \quad \mbox{ holds for } 1 \leq r \leq k.
	  \end{align*}
	 \end{defn}

	\begin{prop}[Koga\cite{Koga}]\label{Koga}
	 \begin{description}
	  \item[(1)]Assume 
		     $u \otimes v = t(a_1,\ldots,a_N) \otimes t(b_1,\ldots,b_N) \in \ovB_{sp}^+ \otimes \ovB_{sp}^+$.
	 Then we have

	 \begin{description}
	  \item[(1A)]For $0 \leq k \leq N-2$, $k \equiv N \mod 2$,
	 \begin{align*}  
	     &u \otimes v \in \ovB(-\xi_k)  
	     \Longleftrightarrow\\
	     &t(a_1,\ldots,a_{N-k};a_{N-k+1},\ldots,a_N|b_1,\ldots,b_k;b_{k+1},\ldots,b_N)
               \mbox{ is semi-standard and }\\
	     &t(a_1,\ldots,a_{N-k-2};a_{N-k-1},\ldots,a_N|b_1,\ldots,b_{k+2};b_{k+3},\ldots,b_N)
	  \mbox{ is not semi-standard,}
	 \end{align*}
	  \item[(1B)]
	 \begin{eqnarray}\nonumber  
	     u \otimes v \in \ovB(-\xi_N)  
	     \Longleftrightarrow
	     t(;a_1,\ldots,a_N|b_1,\ldots,b_N;)
	     \mbox{ is semi-standard. }
	 \end{eqnarray}
	 \end{description}
	 \item[(2)]Assume
		     $u \otimes v = t(a_1,\ldots,a_N) \otimes t(b_1,\ldots,b_N) \in \ovB_{sp}^+ \otimes \ovB_{sp}^-$.
		     Then we have\\
	 \begin{description}
	 \item[(2A)]For $0 \leq k \leq N-3$, $k \equiv N-1 \mod 2$
	 \begin{align*}  
	     &u \otimes v \in \ovB(-\xi_k)  
	     \Longleftrightarrow\\
	     &t(a_1,\ldots,a_{N-k};a_{N-k+1},\ldots,a_N|b_1,\ldots,b_k;b_{k+1},\ldots,b_N)
	  \mbox{ is semi-standard and }\\
	     &t(a_1,\ldots,a_{N-k-2};a_{N-k-1},\ldots,a_N|b_1,\ldots,b_{k+2};b_{k+2},\ldots,b_N)
	  \mbox{ is not semi-standard,}
	 \end{align*}
	 
	 \item[(2B)]
	 \begin{eqnarray}\nonumber  
	     u \otimes v \in \ovB(-\xi_{N-1})  
	     \Longleftrightarrow
	     t(a_1;a_2,\ldots,a_N|b_1,\ldots,b_{N-1};b_N)
	     \mbox{ is semi-standard. }
	 \end{eqnarray}
	\end{description}
	\end{description}
	\end{prop}

\section{Crystal Base of {\boldmath$V(-\omega_N)$}}

	We describe an analogue of spin representations for $U_q(D(N,1))$
	using the $\udn$-crystals.
	
	\begin{prop}\label{spin}
		The irreducible lowest weight module  $V(-\omega_N)$ with lowest weight
	        $-\omega_N$ has a basis over ${\mathbb{Q}}(q)$
	 
		\begin{equation}\nonumber
		 \begin{split}
		&\{{v}(i_1,\ldots,i_N)_{2n} \>; \> n \in \pos, (i_1,\ldots,i_N)
		\in \ovB_{sp}^+\}\sqcup\\
		&\{{v}(i_1,\ldots,i_N)_{2n+1} \>; \> n \in \pos, (i_1,\ldots,i_N)
		\in \ovB_{sp}^-\}
		 \end{split}
		\end{equation}
		with the lowest weight vector ${v}(+,\ldots,+)_0$
	        such that the actions of $\sigma$ and $e_i$ read as follows;
	        \begin{eqnarray}
		 \sigma {v}(+,\ldots,+)_0 = {v}(+,\ldots,+)_0, 
		\end{eqnarray}
	        \begin{eqnarray}\label{eiact}
		  {e_i}({v}(i_1,\ldots,i_N)_k) =\begin{cases}
		     {v}(i_1',\ldots,i_N')_k &\mbox{ if }\widetilde{e_i}(i_1,\ldots,i_N)=
		                                              (i_1',\ldots,i_N') \neq 0 \text{ in }\ovB_{sp}^\pm,\\
		     0                             &\mbox{ otherwise, }
		    \end{cases}
		\end{eqnarray}
		for $1 \leq i \leq N$, and
	        \begin{equation}\label{e_0action}
		  \hspace*{-20mm}{e_0}({v}(i_1,i_2,\ldots,i_N)_k) =
		    \begin{cases}
		     q^{-k}v(+,i_2,\ldots,i_N)_{k+1} &\mbox{ if }i_1=-,\\
		     0                                 &\mbox{ otherwise. }
		    \end{cases}
		\end{equation}
	 \begin{proof}
	  Let ${v}(+,\ldots,+)_0$ be the lowest weight vector of
	  $V(-\omega_N)$.
	  \\

	  Claim 1: $V(-\omega_N)$ is infinite dimensional.

	  As $U_q'(D(N,1)_0) \supset U_q(C_1)$-module, the weight of ${v}(+,\ldots,+)_0$
	  is $1$. Hence $V(-\omega_N)$ is infinite dimensional.
          \\

	  Because $-\omg_N$ is $-\Lam_N$ as a weight of $\udn$,
	  $v(+,\ldots,+)_0$ is the lowest weight vector of the spin representation.
	  For $(i_1',\ldots,i_N') \in \ovB_{sp}^+$, we define
	  ${v}(i_1',\ldots,i_N')_0$ by
	  \begin{equation}\label{claim2}
	   {v}(i_1',\ldots,i_N')_0 = {e_i}({v}(i_1,\ldots,i_N)_0) \quad 
	    \mbox{ where }\widetilde{e_i}(i_1,\ldots,i_N)=(i_1',\ldots,i_N'). 
	  \end{equation}
	  In $\ovV(-\Lam_N)$, $\tei = e_i$ on each weight vector for $1 \leq i \leq N$.
	  Moreover, if $\tei \widetilde{e_j}=\widetilde{e_j} \tei$ for
	  $1 \leq i \neq j \leq N$ in $\ovB_{sp}^+$, $i$-th node and $j$-th
	  node are not connected in the Dynkin diagram.
	  Hence $e_ie_j = e_je_i$.
	  These mean that (\ref{claim2}) is well-defined.
	  (\ref{eiact}) holds for $k=0$ by the definition.
	  \\

%

	  Claim 2: $e_0(v(+,\ldots,+)_0) = 0$

	  Because $f_i(e_0(v(+,\ldots,+)_0)) = 0$ for
	  $0 \leq i \leq N$ by (\ref{efcom}), $e_0(v(+,\ldots,+)_0)$ is a singular vector if it is not
	  $0$. This contradicts the irreducibility of $V(-\omega_N)$. 
	  \\

	  Claim 3: $e_0(v(+,i_2,\ldots,i_N)_0) = 0$

	  Let $v(+,i_2,\ldots,i_N)_0 = e_{l_1} \cdots e_{l_p}(v(+,\ldots,+)_0)$.
	  By (\ref{claim2}), 
	  $(+,i_2,\ldots,i_N) = \widetilde{e_{l_1}} \cdots \widetilde{e_{l_p}}(+,\ldots,+)$ in $\ovB_{sp}^+$.
	  Because $\widetilde{e_1}$ changes $(+,-,\ldots)$ into
	  $(-,+,\ldots)$ in $\ovB_{sp}^+$, it follows that $l_1,\ldots,l_p \in \{2,\ldots,N\}$.
	  By (\ref{commutation}) and Claim 2, we have 
	  \begin{equation}\label{used}
	   \begin{split}
	  e_0(v(+,i_2,\ldots,i_N)_0) &= e_0e_{l_1} \cdots e_{l_p}(v(+,\ldots,+)_0)\\
	                                     &= e_{l_1} \cdots e_{l_p}e_0(v(+,\ldots,+)_0)\\
	                                     &= 0.
	   \end{split}
	  \end{equation}

	  Claim 4: $e_0(v(-,+,\ldots,+,-)_0) \neq 0$

	  By Claim 1 and Claim 3, there exists
	  $v(-,i_2,\ldots,i_N)_0$ such that
	  $e_0(v(-,i_2,\ldots,i_N)_0) \neq 0$.
	  Because $-$ appears an even number of times in 
	  $(-,i_2,\ldots,i_N) \in \ovB_{sp}^+$, 
	  we may assume 
	  \begin{equation}\nonumber
	   \begin{split}
	    (-,i_2,\ldots,i_N) = \widetilde{e_{l_1}} \cdots \widetilde{e_{l_r}}(-,+,\ldots,+,-)\\
	   \mbox{ with }\quad l_1, \ldots ,l_r \in \{2,\ldots,N\}.
	   \end{split}
	  \end{equation}
	  This implies
	  \begin{equation}\nonumber
	   \begin{split}
	   v(-,i_2,\ldots,i_N)_0 = e_{l_1} \cdots e_{l_r}v(-,+,\ldots,+,-)_0\\
	   \mbox{ with }\quad l_1, \ldots ,l_r \in \{2,\ldots,N\}
	   \end{split}
	  \end{equation}
	  by (\ref{claim2}). Hence,
	  \begin{equation}
	   \begin{split}\nonumber
	    0 \neq e_0(v(-,i_2,\ldots,i_N)_0) &= e_0 e_{l_1} \cdots e_{l_r}v(-,+,\ldots,+,-)_0\\
	                                              &= e_{l_1} \cdots e_{l_r}e_0v(-,+,\ldots,+,-)_0.
	   \end{split}
	  \end{equation}
	  
	  We put
	  \begin{align}\label{casek}
	   v(+,\ldots,+,-)_1 = e_0(v(-,+,\ldots,+,-)_0).
	  \end{align}
	  Since $f_i(v(-,+,\ldots,+,-)_0) = 0$ holds for $2 \leq i \leq N$, 
	  we have 

	   \begin{align*}
	    f_iv(+,\ldots,+,-)_1 &= f_ie_0v(-,+,\ldots,+,-)_0 \\
	                         &= 0\quad \mbox{ for } 1 \leq i \leq N,
	   \end{align*}
	  where we used (\ref{used}) for $i=1$.
	  This implies that
	  $v(+,+,\ldots,+,-)_1$ is the lowest weight vector with
	  lowest weight $-\Lambda_{N-1}$ as $\udn$-module. As in the
	  case of $k=0$,  we define
	  ${v}(i_1',\ldots,i_N')_1$ for $(i_1',\ldots,i_N') \in \ovB_{sp}^-$ by
	  \begin{equation}\nonumber
	   {v}(i_1',\ldots,i_N')_1 = {e_i}({v}(i_1,\ldots,i_N)_1), \quad 
	    \mbox{ where }\widetilde{e_i}(i_1,\ldots,i_N)=(i_1',\ldots,i_N'). 
	  \end{equation}
	  Then (\ref{eiact}) holds for $k=1$.
	  \\

	  Claim 5: \begin{equation}\nonumber
	   e_0(v(-,i_2,\ldots,i_N)_0) = v(+,i_2,\ldots,i_N)_1 \quad \mbox{ for any }(-,i_2,\ldots,i_N) \in \ovB_{sp}^+
	  \end{equation}
	  
	  By the definition of $v(+,i_2,\ldots,i_N)_1$, we have
	  \begin{align}\nonumber
	   v(+,i_2,\ldots,i_N)_1 &= e_{l_1} \cdots
	                            e_{l_r}e_0v(-,+,\ldots,+,-)_0 \quad
	   \mbox{ with }
	                             l_1,\ldots,l_r \in \{2,\ldots,N\}\\\nonumber
	                         &= e_0e_{l_1} \cdots e_{l_r}v(-,+,\ldots,+,-)_0\\\nonumber
	                         &= e_0(v(-,i_2,\ldots,i_N)_0).
	  \end{align}

	  In the cases of $k \geq 1$, we put
	  \begin{equation}\nonumber
	   v(+,\ldots,+,-)_{k+1} = q^{-k}e_0(v(-,+,\ldots,+,-)_k)
	  \end{equation}
	  in place of (\ref{casek}). Then the rest of the proof is similar.
	 \end{proof}
	\end{prop}

	\begin{rem}\label{coeff}
	 Suppose $wt(v(i_1,\ldots,i_N)_k)=n_0\omega_0 - \sum_{i=1}^{N}n_i\omega_i$. Then
	 \begin{equation}\nonumber
	  n_0 = \begin{cases}
		 k    &\mbox{ if }i_1 = +,\\
		 k+1  &\mbox{ if }i_1 = -.
		\end{cases}
	 \end{equation}
	 This is because among $e_i$'s only $e_1$ changes the value of $n_0$.
	\end{rem}


	\begin{prop}\label{superspin}
	 The irreducible lowest weight module $V(-{\omega}_N)$ 
	 has a polarizable crystal base $(L,B)$ given as follows.
		\begin{eqnarray}\label{phi}
		L = 	\bigoplus_{\substack{(i_1,\ldots,i_N) \in \ovB_{sp}^+\\
				n \in \pos }}
				Av(i_1,\ldots,i_N)_{2n}\oplus
			\bigoplus_{\substack{(i_1,\ldots,i_N) \in \ovB_{sp}^-\\
				n \in \pos }}
				Av(i_1,\ldots,i_N)_{2n+1}
		\end{eqnarray}
		\begin{eqnarray}
		 \begin{split}
		B =& \{\pm v(i_1,\ldots,i_N)_{2n} \mod qL ; (i_1,\ldots,i_N) \in \ovB_{sp}^+, n \in \pos\} \sqcup\\
		   & \{\pm v(i_1,\ldots,i_N)_{2n+1} \mod qL ; (i_1,\ldots,i_N) \in \ovB_{sp}^-, n \in \pos\}
		 \end{split}
		\end{eqnarray}
	 
	 The Kashiwara operators $e_i$ act on $B$ as (we omit $\mod qL$)
	 \begin{eqnarray}
	  \label{kashiwaraonspin}
	    \widetilde{e_i}v(i_1,\ldots,i_N)_k&=&
		    \begin{cases}
		     v(i_1',\ldots,i_N')_k &\mbox{ if }\widetilde{e_i}(i_1,\ldots,i_N)=
		                                              (i_1',\ldots,i_N') \neq 0 \text{ in }\ovB_{sp}^\pm,\\
		     0                                  &\mbox{ otherwise, }
		    \end{cases}\\\nonumber
		\mbox{ for $1 \leq i \leq N$ and,}\\\label{kashiwaraonspin0}
	    \widetilde{e_0}v(i_1,i_2,\ldots,i_N)_k&=&
		    \begin{cases}
		     v(+,i_2,\ldots,i_N)_{k+1} &\mbox{ if }i_1=-,\\
		     0                                  &\mbox{ otherwise. }
		    \end{cases}
	 \end{eqnarray}

	 \begin{proof}
	   First, we show that $L$ is a crystal lattice.

           It suffices to show (L3). We have only to show it in the case
	  of $i=0$ because the remaining cases are the same as in $\udn$. 
	  Since $\langle h_0,wt(v(-,i_2,\ldots,i_N)_k) \rangle = k+1$, 
	  \begin{eqnarray}\label{e_0tilde}
	   \begin{split}
	   \widetilde{e_0}v(-,i_2,\ldots,i_N)_k &= q_0^{-1}t_0e_0v(-,i_2,\ldots,i_N)_k\\ 
	                                    &= q_0^{-1}q_0^{k+1}q^{-k}v(+,i_2,\ldots,i_N)_{k+1}\\
	                                    &= v(+,i_2,\ldots,i_N)_{k+1}
	   \end{split}
	  \end{eqnarray}
	  and
	  \begin{eqnarray}\label{f_0tilde}
	   \begin{split}
	   \widetilde{f_0}v(+,i_2,\ldots,i_N)_{k+1} &= f_0v(+,i_2,\ldots,i_N)_{k+1}\\ 
	                                        &= q^kf_0e_0v(-,i_2,\ldots,i_N)_{k}\\
	                                        &= \frac{q^{2k+2}-1}{q^2-1}v(-,i_2,\ldots,i_N)_{k} \in L.
	   \end{split}
	  \end{eqnarray}

	   Next, we show that $(L,B)$ is a crystal base. (B1) and (B2)
	  follow from the definition of $B$. (B3) and (B4)
	  follow from (\ref{e_0tilde}) and (\ref{f_0tilde}).

	   Finally the following symmetric bilinear form on
	  $V(-\omega_N)$ is a polarization.
	 
	  \begin{align}\nonumber
	   \left(v(i_1,\ldots,i_N)_0,v(i_1,\ldots,i_N)_0\right) &= 1,\\\nonumber
	   \left(v(i_1,\ldots,i_N)_k,v(i_1,\ldots,i_N)_k\right) &=
	   \prod_{j=1}^{k}\frac{q^{2j}-1}{q^2-1} \quad \mbox{ if } k \geq 1,\\\nonumber
	                                                        &=0
	   \mbox{\quad otherwise.}
	  \end{align}
	  Hence, $(L,B)$ is a polarizable crystal base.
	 \end{proof}
		
	\end{prop}

	 Figure 1 is the crystal graph of $B(-\omega_4)$ of $D(4,1)$ 
	(we omit $v$ and $\mod qL$). 
	The 8 nodes connected with each other by
	$i$-arrow ($1 \leq i \leq 4$) form the crystal graph of 
	$\ovB_{sp}^{\pm}$ of $U_q(D(4))$. 
	The rightmost $0$-arrow for example changes
	the first signature from $+$ into $-$ and 
	the integer from $0$ into $1$.
	The coefficient of $\ome0$ increases $1$ each time we cross $1$-arrow
	from the lower right to the upper left.

	\vspace{2mm}
	\hspace*{-2mm}
	\begin{center}
\unitlength 0.1in
\begin{picture}( 61.1220, 31.4961)(  3.8386,-39.3701)
%
\special{pn 20}%
\special{sh 1}%
\special{ar 6497 3544 10 10 0  6.28318530717959E+0000}%
\special{sh 1}%
\special{ar 5906 3544 10 10 0  6.28318530717959E+0000}%
\special{sh 1}%
\special{ar 5315 3544 10 10 0  6.28318530717959E+0000}%
\special{sh 1}%
\special{ar 4725 3150 10 10 0  6.28318530717959E+0000}%
\special{sh 1}%
\special{ar 4725 3938 10 10 0  6.28318530717959E+0000}%
\special{sh 1}%
\special{ar 4134 3544 10 10 0  6.28318530717959E+0000}%
\special{sh 1}%
\special{ar 3544 3544 10 10 0  6.28318530717959E+0000}%
\special{sh 1}%
\special{ar 2953 3544 10 10 0  6.28318530717959E+0000}%
\special{sh 1}%
\special{ar 2953 3544 10 10 0  6.28318530717959E+0000}%
%
\special{pn 13}%
\special{pa 5965 3544}%
\special{pa 6447 3544}%
\special{fp}%
\special{sh 1}%
\special{pa 6447 3544}%
\special{pa 6381 3524}%
\special{pa 6395 3544}%
\special{pa 6381 3563}%
\special{pa 6447 3544}%
\special{fp}%
\special{pa 5365 3544}%
\special{pa 5837 3544}%
\special{fp}%
\special{sh 1}%
\special{pa 5837 3544}%
\special{pa 5771 3524}%
\special{pa 5785 3544}%
\special{pa 5771 3563}%
\special{pa 5837 3544}%
\special{fp}%
\special{pa 4764 3170}%
\special{pa 5286 3544}%
\special{fp}%
\special{sh 1}%
\special{pa 5286 3544}%
\special{pa 5244 3490}%
\special{pa 5244 3513}%
\special{pa 5221 3521}%
\special{pa 5286 3544}%
\special{fp}%
\special{pa 4754 3908}%
\special{pa 5286 3563}%
\special{fp}%
\special{sh 1}%
\special{pa 5286 3563}%
\special{pa 5220 3582}%
\special{pa 5242 3592}%
\special{pa 5242 3616}%
\special{pa 5286 3563}%
\special{fp}%
\special{pa 4174 3524}%
\special{pa 4705 3170}%
\special{fp}%
\special{sh 1}%
\special{pa 4705 3170}%
\special{pa 4639 3189}%
\special{pa 4662 3199}%
\special{pa 4662 3223}%
\special{pa 4705 3170}%
\special{fp}%
\special{pa 4174 3573}%
\special{pa 4676 3918}%
\special{fp}%
\special{sh 1}%
\special{pa 4676 3918}%
\special{pa 4632 3865}%
\special{pa 4632 3888}%
\special{pa 4611 3897}%
\special{pa 4676 3918}%
\special{fp}%
\special{pa 3593 3544}%
\special{pa 4085 3544}%
\special{fp}%
\special{sh 1}%
\special{pa 4085 3544}%
\special{pa 4019 3524}%
\special{pa 4033 3544}%
\special{pa 4019 3563}%
\special{pa 4085 3544}%
\special{fp}%
\special{pa 2983 3544}%
\special{pa 3495 3544}%
\special{fp}%
\special{sh 1}%
\special{pa 3495 3544}%
\special{pa 3429 3524}%
\special{pa 3442 3544}%
\special{pa 3429 3563}%
\special{pa 3495 3544}%
\special{fp}%
%
\special{pn 20}%
\special{sh 1}%
\special{ar 4725 2560 10 10 0  6.28318530717959E+0000}%
\special{sh 1}%
\special{ar 4134 2560 10 10 0  6.28318530717959E+0000}%
\special{sh 1}%
\special{ar 3544 2560 10 10 0  6.28318530717959E+0000}%
\special{sh 1}%
\special{ar 2953 2166 10 10 0  6.28318530717959E+0000}%
\special{sh 1}%
\special{ar 2953 2953 10 10 0  6.28318530717959E+0000}%
\special{sh 1}%
\special{ar 2363 2560 10 10 0  6.28318530717959E+0000}%
\special{sh 1}%
\special{ar 1772 2560 10 10 0  6.28318530717959E+0000}%
\special{sh 1}%
\special{ar 1182 2560 10 10 0  6.28318530717959E+0000}%
\special{sh 1}%
\special{ar 1182 2560 10 10 0  6.28318530717959E+0000}%
%
\special{pn 13}%
\special{pa 4193 2560}%
\special{pa 4676 2560}%
\special{fp}%
\special{sh 1}%
\special{pa 4676 2560}%
\special{pa 4610 2540}%
\special{pa 4624 2560}%
\special{pa 4610 2579}%
\special{pa 4676 2560}%
\special{fp}%
\special{pa 3593 2560}%
\special{pa 4065 2560}%
\special{fp}%
\special{sh 1}%
\special{pa 4065 2560}%
\special{pa 4000 2540}%
\special{pa 4013 2560}%
\special{pa 4000 2579}%
\special{pa 4065 2560}%
\special{fp}%
\special{pa 2993 2186}%
\special{pa 3514 2560}%
\special{fp}%
\special{sh 1}%
\special{pa 3514 2560}%
\special{pa 3472 2505}%
\special{pa 3472 2529}%
\special{pa 3449 2537}%
\special{pa 3514 2560}%
\special{fp}%
\special{pa 2983 2924}%
\special{pa 3514 2579}%
\special{fp}%
\special{sh 1}%
\special{pa 3514 2579}%
\special{pa 3448 2598}%
\special{pa 3470 2608}%
\special{pa 3470 2631}%
\special{pa 3514 2579}%
\special{fp}%
\special{pa 2402 2540}%
\special{pa 2934 2186}%
\special{fp}%
\special{sh 1}%
\special{pa 2934 2186}%
\special{pa 2868 2205}%
\special{pa 2890 2215}%
\special{pa 2890 2239}%
\special{pa 2934 2186}%
\special{fp}%
\special{pa 2402 2589}%
\special{pa 2904 2934}%
\special{fp}%
\special{sh 1}%
\special{pa 2904 2934}%
\special{pa 2861 2880}%
\special{pa 2861 2904}%
\special{pa 2839 2913}%
\special{pa 2904 2934}%
\special{fp}%
\special{pa 1821 2560}%
\special{pa 2313 2560}%
\special{fp}%
\special{sh 1}%
\special{pa 2313 2560}%
\special{pa 2248 2540}%
\special{pa 2261 2560}%
\special{pa 2248 2579}%
\special{pa 2313 2560}%
\special{fp}%
\special{pa 1211 2560}%
\special{pa 1723 2560}%
\special{fp}%
\special{sh 1}%
\special{pa 1723 2560}%
\special{pa 1657 2540}%
\special{pa 1671 2560}%
\special{pa 1657 2579}%
\special{pa 1723 2560}%
\special{fp}%
%
\special{pn 20}%
\special{sh 1}%
\special{ar 2953 1575 10 10 0  6.28318530717959E+0000}%
\special{sh 1}%
\special{ar 2363 1575 10 10 0  6.28318530717959E+0000}%
\special{sh 1}%
\special{ar 1772 1575 10 10 0  6.28318530717959E+0000}%
\special{sh 1}%
\special{ar 1182 1969 10 10 0  6.28318530717959E+0000}%
\special{sh 1}%
\special{ar 1182 1182 10 10 0  6.28318530717959E+0000}%
\special{sh 1}%
\special{ar 591 1575 10 10 0  6.28318530717959E+0000}%
\special{sh 1}%
\special{ar 1772 1565 10 10 0  6.28318530717959E+0000}%
%
\special{pn 13}%
\special{pa 2402 1575}%
\special{pa 2904 1575}%
\special{fp}%
\special{sh 1}%
\special{pa 2904 1575}%
\special{pa 2838 1556}%
\special{pa 2852 1575}%
\special{pa 2838 1595}%
\special{pa 2904 1575}%
\special{fp}%
\special{pa 1802 1575}%
\special{pa 2313 1575}%
\special{fp}%
\special{sh 1}%
\special{pa 2313 1575}%
\special{pa 2248 1556}%
\special{pa 2261 1575}%
\special{pa 2248 1595}%
\special{pa 2313 1575}%
\special{fp}%
\special{pa 1221 1211}%
\special{pa 1733 1536}%
\special{fp}%
\special{sh 1}%
\special{pa 1733 1536}%
\special{pa 1688 1484}%
\special{pa 1688 1507}%
\special{pa 1667 1517}%
\special{pa 1733 1536}%
\special{fp}%
\special{pa 1221 1939}%
\special{pa 1743 1585}%
\special{fp}%
\special{sh 1}%
\special{pa 1743 1585}%
\special{pa 1678 1606}%
\special{pa 1699 1615}%
\special{pa 1699 1638}%
\special{pa 1743 1585}%
\special{fp}%
%
\special{pn 13}%
\special{pa 621 1546}%
\special{pa 1142 1191}%
\special{fp}%
\special{sh 1}%
\special{pa 1142 1191}%
\special{pa 1077 1212}%
\special{pa 1099 1221}%
\special{pa 1099 1245}%
\special{pa 1142 1191}%
\special{fp}%
\special{pa 611 1605}%
\special{pa 1142 1969}%
\special{fp}%
\special{sh 1}%
\special{pa 1142 1969}%
\special{pa 1099 1916}%
\special{pa 1099 1939}%
\special{pa 1077 1948}%
\special{pa 1142 1969}%
\special{fp}%
%
\special{pn 8}%
\special{pa 4725 2560}%
\special{pa 4725 3150}%
\special{dt 0.045}%
\special{sh 1}%
\special{pa 4725 3150}%
\special{pa 4745 3084}%
\special{pa 4725 3098}%
\special{pa 4705 3084}%
\special{pa 4725 3150}%
\special{fp}%
\special{pa 4134 2560}%
\special{pa 4134 3544}%
\special{dt 0.045}%
\special{sh 1}%
\special{pa 4134 3544}%
\special{pa 4154 3478}%
\special{pa 4134 3492}%
\special{pa 4115 3478}%
\special{pa 4134 3544}%
\special{fp}%
\special{pa 3544 2560}%
\special{pa 3544 3544}%
\special{dt 0.045}%
\special{sh 1}%
\special{pa 3544 3544}%
\special{pa 3563 3478}%
\special{pa 3544 3492}%
\special{pa 3524 3478}%
\special{pa 3544 3544}%
\special{fp}%
\special{pa 2953 2953}%
\special{pa 2953 3544}%
\special{dt 0.045}%
\special{sh 1}%
\special{pa 2953 3544}%
\special{pa 2973 3478}%
\special{pa 2953 3492}%
\special{pa 2934 3478}%
\special{pa 2953 3544}%
\special{fp}%
\special{pa 2953 1575}%
\special{pa 2953 2166}%
\special{dt 0.045}%
\special{sh 1}%
\special{pa 2953 2166}%
\special{pa 2973 2100}%
\special{pa 2953 2114}%
\special{pa 2934 2100}%
\special{pa 2953 2166}%
\special{fp}%
\special{pa 2363 1575}%
\special{pa 2363 2560}%
\special{dt 0.045}%
\special{sh 1}%
\special{pa 2363 2560}%
\special{pa 2382 2494}%
\special{pa 2363 2507}%
\special{pa 2343 2494}%
\special{pa 2363 2560}%
\special{fp}%
\special{pa 1772 1575}%
\special{pa 1772 2560}%
\special{dt 0.045}%
\special{sh 1}%
\special{pa 1772 2560}%
\special{pa 1792 2494}%
\special{pa 1772 2507}%
\special{pa 1752 2494}%
\special{pa 1772 2560}%
\special{fp}%
\special{pa 1182 1969}%
\special{pa 1182 2560}%
\special{dt 0.045}%
\special{sh 1}%
\special{pa 1182 2560}%
\special{pa 1201 2494}%
\special{pa 1182 2507}%
\special{pa 1162 2494}%
\special{pa 1182 2560}%
\special{fp}%
\special{pa 1182 788}%
\special{pa 1182 1182}%
\special{dt 0.045}%
\special{sh 1}%
\special{pa 1182 1182}%
\special{pa 1201 1116}%
\special{pa 1182 1129}%
\special{pa 1162 1116}%
\special{pa 1182 1182}%
\special{fp}%
\special{pa 591 1182}%
\special{pa 591 1575}%
\special{dt 0.045}%
\special{sh 1}%
\special{pa 591 1575}%
\special{pa 611 1509}%
\special{pa 591 1523}%
\special{pa 571 1509}%
\special{pa 591 1575}%
\special{fp}%
\put(59.0551,-34.9409){\makebox(0,0)[lb]{$(+,+,+,+)_0$}}%
\put(48.2283,-25.0984){\makebox(0,0)[lb]{$(+,+,+,-)_1$}}%
\put(30.5118,-15.1575){\makebox(0,0)[lb]{$(+,+,+,+)_2$}}%
\put(48.2283,-31.2008){\makebox(0,0)[lb]{$(-,+,+,-)_0$}}%
\put(30.5118,-21.1614){\makebox(0,0)[lb]{$(-,+,+,+)_1$}}%
\put(26.5748,-23.6220){\makebox(0,0)[rb]{$4$}}%
\put(32.4803,-23.6220){\makebox(0,0)[lb]{$1$}}%
\put(32.4803,-27.5591){\makebox(0,0)[lt]{$4$}}%
\put(26.5748,-27.5591){\makebox(0,0)[rt]{$1$}}%
\put(14.7638,-17.7165){\makebox(0,0)[lt]{$3$}}%
\put(14.7638,-13.7795){\makebox(0,0)[lb]{$1$}}%
\put(50.1969,-33.4646){\makebox(0,0)[lb]{$1$}}%
\put(50.1969,-37.4016){\makebox(0,0)[lt]{$3$}}%
\put(44.2913,-37.4016){\makebox(0,0)[rt]{$1$}}%
\put(44.2913,-33.4646){\makebox(0,0)[rb]{$3$}}%
\put(44.2913,-27.5591){\makebox(0,0)[rb]{$3$}}%
%
\special{pn 13}%
\special{sh 1}%
\special{ar 522 1575 10 10 0  6.28318530717959E+0000}%
\special{sh 1}%
\special{ar 453 1575 10 10 0  6.28318530717959E+0000}%
\special{sh 1}%
\special{ar 384 1575 10 10 0  6.28318530717959E+0000}%
\special{sh 1}%
\special{ar 384 1595 10 10 0  6.28318530717959E+0000}%
\put(38.3858,-27.5591){\makebox(0,0)[rb]{$2$}}%
\put(38.3858,-37.4016){\makebox(0,0)[rb]{$2$}}%
\put(56.1024,-37.4016){\makebox(0,0)[rb]{$2$}}%
\put(62.0079,-37.4016){\makebox(0,0)[rb]{$4$}}%
\put(32.4803,-35.7283){\makebox(0,0)[lt]{$4$}}%
\put(26.5748,-17.7165){\makebox(0,0)[rb]{$4$}}%
\put(14.8622,-25.8858){\makebox(0,0)[lt]{$3$}}%
\put(20.6693,-17.7165){\makebox(0,0)[rb]{$2$}}%
\put(20.6693,-27.5591){\makebox(0,0)[rb]{$2$}}%
\end{picture}%

	\end{center}
	\begin{center}
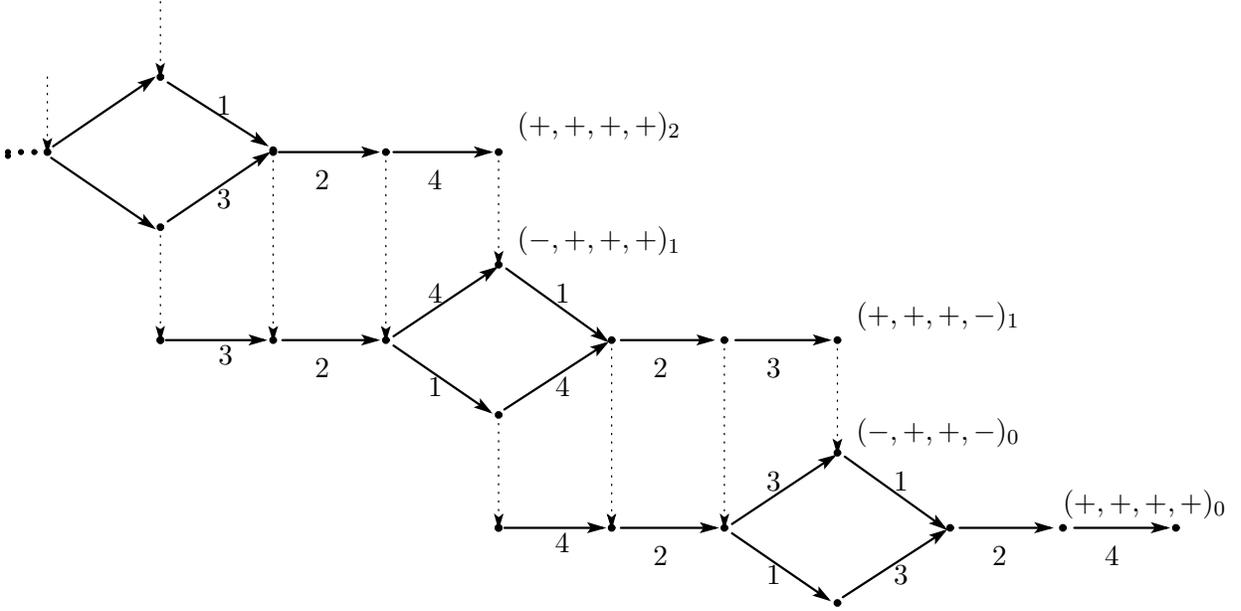
	 Figure 1 The crystal graph of $B(-\omg_4)$
	\end{center}

\section{Crystal Bases for Fundamental Representations}

        It is natural to ask which representations admit crystal base.
	Next theorem, which is one of our
        main results, is the answer.
	Our tool is the decomposition of $\udn$-crystals.
	
	\begin{thm}\label{exist}
	 The irreducible lowest weight module $V(\lambda)$ with
	 the lowest weight
	 \begin{equation}
		\lambda = n_0{\omega}_0 - \sum_{i=1}^{N}n_i{\omega}_i, \quad n_i \in \pos \quad\mbox{ for }
		 0 \leq i \leq N
	 \end{equation} 
        admits a polarizable crystal base.

	\begin{proof}
	 We prove this when $N$ is even. 
	 The proof for odd $N$ is similar.

	 We show that 
	 \begin{eqnarray}\label{spintensor}
	  \begin{split}
	   &LW(B(-\omega_N) \otimes B(-\omega_N)) =\\
	   &\left\{
	   \begin{array}{lll}	   
	    v(+,\ldots,+)_0 \otimes v(-,\ldots,-)_{2k} &                                &k \in \pos,\\
	    v(+,\ldots,+)_0 \otimes v(\overbrace{+,\ldots,+}^{2i},-,\ldots,-)_0& \mbox{ ; } &1 \leq i \leq \frac{N-2}{2}\\
	    v(+,\ldots,+)_0 \otimes v(+,\ldots,+)_0&
	   \end{array}
	   \right\}.\\
	  \end{split}
	 \end{eqnarray}
	 The element in the left hand side of (\ref{spintensor}) is $v(+,\ldots,+)_0 \otimes b$
	 for some $b \in B(-\omega_N)$ by Lemma
	 \ref{lwv}. By (\ref{kashiwaraonspin}) and (\ref{clhwv}) in
	 Proposition \ref{spinspin},
	 \begin{equation}\nonumber
	  \begin{split}
	                         &\quad \tfi\left(v(+,\ldots,+)_0 \otimes b\right) = 0 \quad \mbox{for }1 \leq i \leq N\\
	   \Longleftrightarrow   &\quad b=v(-,\ldots,-)_{2k}  \mbox{ or } v(\overbrace{+,\ldots,+}^{j},-,\ldots,-)_l
	                         \quad\mbox{for some }k,l \in \pos,1 \leq j \leq N.
	  \end{split}
	 \end{equation}

	 If $l\geq1$, and $1 \leq j \leq N$,
	 \begin{equation}\nonumber
	  \begin{split}
	 \tf0(v(+,\ldots,+)_0 \otimes v(\overbrace{+,\ldots,+}^{j},-,\ldots,-)_l) 
	 &=v(+,\ldots,+)_0 \otimes \tf0(v(\overbrace{+,\ldots,+}^{j},-,\ldots,-)_l)\\
 	 &=v(+,\ldots,+)_0 \otimes v(-,\overbrace{+\ldots,+}^{j-1},-,\ldots,-)_{l-1}\\
	 &\neq 0.
	 \end{split}
	 \end{equation}

	 Hence $l=0$ for this case and (\ref{spintensor}) follows.

	 As a consequence, we have
	 \begin{eqnarray}\nonumber
	  B(-\omega_N) \otimes B(-\omega_N) &= B(-2\omega_N) \oplus \bigoplus\limits_{j=1}^{\frac{N-2}{2}}B(-\omega_{2j})
	   \oplus \bigoplus\limits_{k \in \mathbb{Z}_{\geq 0}}B((2k+1)\omega_0).
	 \end{eqnarray}
	 Similarly, we have
	 \begin{eqnarray}\nonumber
	  B(-\omega_N) \otimes B(-\omega_{N-1}) &= B(-\omega_N - \omega_{N-1}) \oplus 
	   \bigoplus\limits_{j=0}^{\frac{N-4}{2}}B(-\omega_{2j +1})
	   \oplus \bigoplus\limits_{k \in \mathbb{Z}_{\geq 0}}B((2k+2)\omega_0).
	 \end{eqnarray}
	 In particular, there are polarizable crystal bases with lowest weights
	 $\omega_0$, $-\omega_1$,$\ldots$, $-\omg_N$. Together with
	 Proposition \ref{tensor rule}, we obtain the desired statement.
	\end{proof}
	\end{thm}

	\begin{cor}\label{complete}
	 Let $\lam$ and $\lam'$ be
	 \begin{equation}
		\lambda = n_0{\omega}_0 - \sum_{i=1}^{N}n_i{\omega}_i,\quad 
		\lambda' = n_0'{\omega}_0 - \sum_{i=1}^{N}n_i'{\omega}_i, 
		 \quad n_i, n_i' \in \pos \quad\mbox{ for }
		 0 \leq i \leq N.
	 \end{equation} 
	 Then the tensor product $V(\lam') \otimes V(\lam)$ is
	 completely reducible.
	 \begin{proof}
	  This is a direct consequence of Proposition \ref{tensor rule}.
	 \end{proof}
	\end{cor}

\section{Properties of $B(\ome0)$}	

	Next we treat the tensor product of modules 
	with weights as in Corollary \ref{complete}. 
	Especially we are interested in the case when modules are
	typical because of the following lemma.

	\begin{lem}\label{typeasy}
	 Let $\lam$ and $\lam'$ be as in Corollary \ref{complete}, and
	 $u_{\lam'}$ be the lowest weight vector of $B(\lam')$.
	 If $\lam'$ is typical, then 
	 \begin{equation*}
	  u_{\lam'} \otimes u \in LW( B(\lam') \otimes B(\lam))
	   \Longleftrightarrow
	  u_{\lam'} \otimes u \in \overline{LW}( B(\lam') \otimes B(\lam))
	 \end{equation*}
	 holds for $u \in B(\lam)$.
	 
	 \begin{proof}
	  Assume $u_{\lam'} \otimes u \in \overline{LW}( B(\lam') \otimes B(\lam))$.
	  By (\ref{ftilde0}), we have
	  $\tf0(u_{\lam'} \otimes u) = \tf0(u_{\lam'}) \otimes u = 0$.
	  Hence $u_{\lam'} \otimes u \in LW( B(\lam') \otimes B(\lam))$.
	 \end{proof}
	\end{lem}

	We further restrict $B(\lam)$ in Lemma \ref{typeasy} also to be typical
	because the crystal base of typical representations
	have nice properties which will be stated in Proposition \ref{direct}.
	We first study the structure of $B(\omega_0)$ because it plays
	an important role to investigate typical representations. 
	 $B(\omega_0)$ can be realized by the embedding in the proof of Theorem
	\ref{exist}, that is 	
	 $B(\omega_0) \hookrightarrow B(-\omega_N) \otimes B(-\omega_N)$
	when N is even, 
	$B(\omega_0) \hookrightarrow B(-\omega_N)\otimes B(-\omega_{N-1})$
	when N is odd. 

	Roughly speaking, the proof of Theorem \ref{exist} shows that $B(-\omega_i)$ is
	the union of infinitely many
	$\ovB(-\xi_k)$'s, where each $\ovB(-\xi_k)$ is connected with
	others by $\widetilde{e_0}$ and $\widetilde{f_0}$.  
	We describe more precisely this situation that a $\udn$-crystal is contained 
	in a $U_q(D(N,1))$-crystal.
	\begin{defn}
	 We assume the following conditions.\\
	 (1) $\lam = n_0{\omega}_0 - \sum_{i=1}^{N}n_i{\omega}_i$, \quad$n_j \in \pos$ for $0 \leq j \leq N$,\\ 
	 (2) $\Lam = -\sum_{i=1}^{N}l_i\Lam_i$, \quad$l_i \in \pos$ for all $i$,\\
	 (3) $b \in B(\lam)$ satisfies $wt(b) = l_0\ome0 -\sum_{i=1}^{N}l_i\omg_i$,\quad $l_0 \in \pos$,\\
	 (4) $\tfi b = 0$ for  $1 \leq i \leq N$.\\
	 Then we define a $\udn$-crystal in $B(\lam)$ by
	 \begin{eqnarray}\nonumber
	  \ovB(\Lam ; l_0) = \left\{\widetilde{e_{i_1}}\cdots\widetilde{e_{i_p}}(b) ; 
						     1 \leq i_1, \ldots, i_p \leq N, p \geq 0 \right\}-\{0\}.
	 \end{eqnarray}
	\end{defn}

	\begin{ntn}
	 In order to relate a weight of $\udn$ with that of $U_q(D(N,1))$, we
	 fix some notations.
	 Let $n_0$,$\ldots$, $n_N$ be non-negative integers.\\
	 (1)For $\lambda = n_0\omega_0 - \sum_{i=1}^{N}n_i\omega_i$,
	 we define a dominant integral weight of $\udn$ by
	 \begin{equation}\nonumber
	  \lambda_{cl}=-\sum_{i=1}^{N}n_i\Lambda_i.
	 \end{equation}
	 (2)For $\Lam = -\sum_{i=1}^{N}n_i\Lam_i$, 
	 we define a weight of $U_q(D(N,1))$ by
	 \begin{equation}\nonumber
	  \Lam_{su}= -\sum_{i=1}^{N}n_i\omg_i.
	 \end{equation}
	\end{ntn}

	Note that
	\begin{equation}\nonumber
	 wt(\overline{LW}(\ovB(\Lam ; l_0))) = \Lam_{su} + l_0 \ome0 \quad \mbox{as }U_q(D(N,1))\mbox{-crystal}.
	\end{equation}

	We now determine the places where
	$0$-arrows exist in the above tensor products. 

	\begin{defn}
	 We define
	 $$\begin{CD}
	    \ovB(\Lam;l) @>{0}>R> \ovB(\Lam';l')
	   \end{CD}$$
	 by the condition
	 \begin{eqnarray}\nonumber
	  \mbox{ for any } b \otimes b'\in \ovB(\Lambda;l) \subset 
	   B(-\omega_N) \otimes B(-\omega_N) \mbox{ or  }B(-\omega_N) \otimes B(-\omega_{N-1}),\\ \nonumber
	  \widetilde{f_0}(b \otimes b') = \sigma b \otimes {\tf0 b'} \neq 0
	  \Longrightarrow \widetilde{f_0}(b \otimes b') \in \ovB(\Lambda';l'),
	 \end{eqnarray}
	 and
	 $$\begin{CD}
	  \ovB(\Lambda;l) @>{0}>L> \ovB(\Lambda';l')
	 \end{CD}$$
	 by the condition
	 \begin{eqnarray}\nonumber
	  \mbox{ for any } b \otimes b' \in \ovB(\Lambda;l) \subset 
	   B(-\omega_N) \otimes B(-\omega_N) \mbox{ or }B(-\omega_N) \otimes B(-\omega_{N-1}),\\\nonumber 
	   \widetilde{f_0}(b \otimes b') = (\tf0 b) \otimes b' \neq 0
	  \Longrightarrow \widetilde{f_0}(b \otimes b') \in \ovB(\Lambda';l').
	 \end{eqnarray}	
	\end{defn}
	
	\begin{lem}\label{0connection}
	 In $B(-\omega_N) \otimes B(-\omega_N)$ and $B(-\omega_N) \otimes B(-\omega_{N-1})$, we have\\
	 \begin{description}
	 \item[(1)]for $k=2,3,\ldots,N$
	 $$\begin{CD}
	  \ovB(-\xi_k;l) \hspace{0mm}       @>{0}>R> \ovB(-\xi_{k-1};l-1),
	 \end{CD}$$

	 \item[(2)]
	 $$\begin{CD}	  
	    \hspace*{-10mm}
	    \ovB(-\xi_1;l) \hspace{0mm}     @>{0}>R> \ovB(-\xi_0;l),
	 \end{CD}$$

	 \item[(3)]for $k=2,\ldots,N-1$
	$$\begin{CD}
	   \hspace*{-0mm}
	  \ovB(-\xi_{k-1};l) \hspace{0mm}   @>{0}>L> \ovB(-\xi_k;l-1),
	 \end{CD}$$
	
	 \item[(4)]
	$$\begin{CD}
	   \hspace*{-4mm}
	  \ovB(-\xi_0;l) \hspace{0mm}       @>{0}>L> \ovB(-\xi_1;l-2),
	 \end{CD}$$

	 \item[(5)]
	 $$\begin{CD}
	   \hspace*{2mm}
	  \ovB(-\xi_N;l) \hspace{0mm}       @>{0}>L> \ovB(-\xi_{N-1};l-1),
	   \end{CD}$$

	 \item[(6)]
	 $$\begin{CD}
	   \hspace*{1mm}
	   \ovB(-\xi_{N-1};l) @>{0}>L> \ovB(-\xi_N;l-1),
	   \end{CD}$$
 	 $$\begin{CD}
	   \hspace*{1mm}
	   \ovB(-\xi_{N-1};l) @>{0}>L> \ovB({-\xi_N}';l-1).
	   \end{CD}$$
	 \end{description}
	 \begin{proof}
	  We prove this when $N$ is even.
	  We denote $u \otimes w = v(i_1,\ldots,i_N)_p \otimes v(j_1,\ldots,j_N)_{p'} $ by
	  $t(a_1,\ldots,a_N)_p \otimes t(b_1,\ldots,b_N)_{p'}$.
	  We also assume $p + p' = l$.\\
	
	  Let us show (1) with $1 \leq k \leq N-2$.
	  Assume $u \otimes w  \in \ovB(-\xi_k;l)$ 
	  satisfies $\tf0(u \otimes w) = \sigma u \otimes \tf0 w \neq 0$. 
	  Let us denote $\tf0w = t(b_1',\ldots,b_N')_{p_1}$
	  Because $\tf0$ changes
	  $(+,\ldots)_{p'}$ into $(-,\ldots)_{p'-1}$, $p_1=p'-1$, $b_1 = 1$,
	  $b_l'=b_{l+1}$ ($1 \leq l \leq N-1$), and $b_N' = \overline 1$.

	  By Proposition \ref{Koga}(1A), $u \otimes w \in \ovB(-\xi_k;l)$ implies
	  \begin{equation*}
	   \begin{array}{ll}
	    t(a_1,\ldots,a_{N-k};a_{N-k+1},\ldots,a_N|1,b_2,\ldots,b_k;b_{k+1},\ldots,b_N)&
	    \mbox{ is semi-standard and }\\
	    t(a_1,\ldots,a_{N-k-2};a_{N-k-1},\ldots,a_N|1,b_2,\ldots,b_{k+2};b_{k+3},\ldots,b_N)&
	    \mbox{ is not semi-standard.}
	   \end{array}
	  \end{equation*}
	  It follows that
	  \begin{equation*}
	   \begin{array}{ll}
	    t(a_1,\ldots,a_{N-k+1};a_{N-k+2},\ldots,a_N|b_2,\ldots,b_k;b_{k+1},\ldots,b_N,\overline{1})&
	    \mbox{ is semi-standard and}\\
	    t(a_1,\ldots,a_{N-k-1};a_{N-k},\ldots,a_N|b_2,\ldots,b_{k+2};b_{k+3},\ldots,b_N,\overline{1})&
	    \mbox{ is not semi-standard.}
	   \end{array}
	  \end{equation*}
	  This means that $\tf0(u \otimes w) \in \ovB(-\xi_{k-1};l')$ for some $l'$.
	  Because the lowest weight vector of
	  $\ovB(-\xi_{k-1};l')$ is 
	  $(+,\ldots,+)_p \otimes (\overbrace{+,\ldots,+}^{k-1},-,\ldots,-)_{p'-1}$ as $\udn$-crystal,
	  we have $l'=p + (p'-1) = l-1$ by Remark \ref{coeff}.

	  The rest of (1) and (2) are similar.

	  Let us verify (3).
	  Assume $u \otimes w \in \ovB(-\xi_{k-1};l)$ satisfies 
	  $\tf0(u \otimes w) = (\tf0 u) \otimes w \neq 0$, and 
	  $\tf0u = t(a_1',\ldots,a_N')_{p_2}$.
	  As in (1), $p_2=p-1$, $a_1 = 1$,
	  $a_l'=a_{l+1}$ ($1 \leq l \leq N-1$), $a_N' = \overline 1$.

	  By Proposition \ref{Koga}(2A), 
	  \begin{equation*}
	   \begin{array}{ll}
	    t(1,a_2,\ldots,a_{N-k+1};a_{N-k+2},\ldots,a_N|b_1,\ldots,b_{k-1};b_k,\ldots,b_N)&
	    \mbox{ is semi-standard and }\\
	    t(1,a_2,\ldots,a_{N-k-1};a_{N-k},\ldots,a_N|b_1,\ldots,b_{k+1};b_{k+2},\ldots,b_N)&
	    \mbox{ is not semi-standard.}
	   \end{array}
	  \end{equation*}
	  It follows that
	  \begin{equation*}
	   \begin{array}{ll}
	    t(a_2,\ldots,a_{N-k+1};a_{N-k+2},\ldots,a_N,\overline{1}|b_1,\ldots,b_{k-1},b_k;b_{k+1},\ldots,b_N)&
	    \mbox{ is semi-standard and}\\
	    t(a_2,\ldots,a_{N-k-1};a_{N-k},\ldots,a_N,\overline{1}|b_1,\ldots,b_{k+1},b_{k+2};b_{k+3},\ldots,b_N)&
	    \mbox{ is not semi-standard.}
	   \end{array}
	  \end{equation*}
	  This means that $\tf0(u \otimes w) \in \ovB(-\xi_k;l')$.
	  Because the lowest weight vector of
	  $\ovB(-\xi_k;l')$ is 
	  $(+,\ldots,+)_{p-1} \otimes (\overbrace{+,\ldots,+}^k,-,\ldots,-)_{p'}$ as $\udn$-crystal,
	  we have $l'=l-1$.
	  The proofs of (4),(5) and (6) are similar.
	 \end{proof}
	\end{lem}

	\begin{rem}\label{occur}
	 In Lemma \ref{0connection}, (1) and (2) occur only if
	 $p=0$ by (\ref{ftilde0}).
	\end{rem}
	 
	\begin{cor}\label{ome}
	 We have
	 \begin{eqnarray}\nonumber
	  B(\omega_0) \cong \bigoplus_{\nu \in W}\bigoplus_{i=1}^{i_{\nu}}\bigoplus_{n \in \mathbb{Z}_{\geq 0}}
	   \ovB\left(\nu ; z^i(\nu) +2n\right),
	 \end{eqnarray}
	 where 
	 $W = \left\{-\xi_0,\ldots,-\xi_N,{-\xi_N}'\right\}$. 
	 The integers $i_{\nu},z^i(\nu)$ are given as follows.
	 
	  \begin{equation}\nonumber
	   \begin{array}{lllll}
	    \mbox{ for } \nu = -\xi_0,&  &i_{\nu} = 2,&  z^1(\nu) = 1 , &z^2(\nu) = 2N+1\\
	    \mbox{ for } \nu = -\xi_k \quad (1 \leq k \leq N-1),& 
	                       &i_{\nu} = 2,& z^1(\nu) = k , &z^2(\nu) = 2N-k,\\
	    \mbox{ for } \nu = {-\xi_N}' \mbox{ or } -\xi_N,&  &i_{\nu} = 1,&  z^1(\nu) = N.&\\
	   \end{array}
	  \end{equation}

	\begin{proof}
	 $B(\ome0)$ is the connected component of $B(-\omg_N) \otimes B(-\omg_N)$ 
	  containing $\ovB(0;1)$ when $N$ is even,
	 and of $B(-\omega_N) \otimes B(-\omega_{N-1})$ when $N$ is odd. 
	 This corollary follows from Lemma
	 \ref{0connection} and Remark \ref{occur}.
	\end{proof}
	\end{cor}

%
%

\section{Crystal bases of typical representations}

        The aim of this section is to obtain a description as in Corollary
        \ref{ome} for typical representations.
        In this section, we assume
	\begin{equation}\label{lam}
	 \lam = -\sum_{i=1}^{N}n_i\omg_i,
	 \quad n_i \in \pos \mbox{ for }1 \leq i \leq N.
	\end{equation}
	We denote the lowest weight vector of $B(\lam)$ by $u_{\lam}$.

	When we consider typical representations, 
	the coefficient of $\ome0$ does not matter essentially.
        \begin{lem}\label{lemshift}
	 Let $\lambda$ be as in (\ref{lam}) and $n \geq 1$.
	 Then there is a bijection 
	 \begin{eqnarray}\label{shift}
	  \varphi: B(\lambda + n\omega_0) \xrightarrow{\sim}B(\lambda + (n+1)\omega_0) 
	 \end{eqnarray}
	 which commutes with the Kashiwara operators and
	 \begin{eqnarray}\label{wtshift}
	  wt(\varphi(b)) = wt(b) + \omega_0.
	 \end{eqnarray}
	 \begin{proof}
	  Note that $B(\lambda + (n+1)\omega_0)$ is the lowest component
	  of 
	  $B(\lambda + n\omega_0) \otimes B(\omega_0)$. 
	  We define $\varphi$ by $\varphi(b) = b \otimes u_{\ome0}$.
	  Lemma \ref{lwv} implies $b \otimes u_{\ome0} \in B(\lambda + (n+1)\omega_0)$
	  for all $b \in B(\lam +n\ome0)$. 
	  Hence it suffices to show that $\varphi$
	  is surjective and that $\varphi$ commutes with $\widetilde{e_i}$'s and
	  $\widetilde{f_i}$'s. 
	  To show the surjectivity, we show that for all $i$
	  \begin{equation}\label{stable}
	   \begin{split}
	   \tei\left(\varphi\left(B(\lambda + n\omega_0)\right)\right) \subset 
	    \varphi\left(B(\lambda + n\omega_0)\right),\\
	   \tfi\left(\varphi\left(B(\lambda + n\omega_0)\right)\right) \subset 
	    \varphi\left(B(\lambda + n\omega_0)\right)
	   \end{split}
	  \end{equation}
	  instead. By Proposition \ref{0connection}, $\widetilde{e_i}u_{\ome0}=0$ for
	  $1 \leq i \leq N$ and by (\ref{etildei}), it follows that
	  $\widetilde{e_i}(b \otimes u_{\ome0})=(\widetilde{e_i}b) \otimes u_{\ome0}$
	  for $1 \leq i \leq N$. 
	  Since $\langle h_0, wt(u_{\ome0}) \rangle > 0$, we have
	  $\widetilde{e_0}(b \otimes u_{\ome0})=(\widetilde{e_0}b)\otimes{u_{\ome0}}$. 
	  The case of $\tfi$'s is similar. Hence (\ref{stable}) holds   
	  and $\varphi$ commutes with $\widetilde{e_i}$'s and $\widetilde{f_i}$'s.
    	 \end{proof}
	\end{lem}

        In view of the above lemma, we have only to consider typical
	representations 
	$B(\lambda + \omega_0)$ with
        $\lambda$ as in (\ref{lam}).
	$B(\lambda + \omega_0)$ is the connected
        component of $B(\omega_0) \otimes B(\lambda)$ containing
        $u_{\omega_0} \otimes u_{\lambda}$. Note that $u_{\lambda}$ 
	is contained in $\ovB(\lam_{cl};0)$.

	The next proposition is one of the favorable properties of
	typical representations.
	\begin{prop}\label{direct}
	 Let $\lambda$ be as in (\ref{lam}). Then we have
	 \begin{eqnarray}\label{typtab}
	  B(\lambda + \omega_0) \cong \left\{ b \otimes u ;  b \in B(\ome0),u \in \ovB(\lambda_{cl};0)
				   \right\}.
	 \end{eqnarray}
	 
	 \begin{proof}
	  Let $J$ be the right hand side of (\ref{typtab}). First we claim that
	  $J$ is stable under $\tei$'s and $\tfi$'s.
	  Since $\ovB(\lambda_{cl};0) \subset B(\lambda)$ is stable under $\tei$,$\tfi$
	  ($1 \leq i \leq N$), $J$ is stable under $\tei$,$\tfi$
	  ($1 \leq i \leq N$). 
	  Because
	  \begin{align}
	   \te0(b \otimes u)&=(\te0 b) \otimes u,\\\label{0apl}
	   \tf0(b \otimes u)&=(\tf0 b) \otimes u
	  \end{align}
	  holds for $b \otimes u \in J$,
	  $J$ is stable under all $\tei$ and $\tfi$.

	  Next we show that $LW(J) = \{u_{\omega_0} \otimes u_{\lambda}\}$. 
	  Assume $LW(J) \ni b \otimes u$ satisfies
	  $b \neq u_{\omega_0}$ and $u \neq u_{\lambda}$. 
	  There exists $0 \leq i \leq N$ such
	  that $\tfi b \neq 0$. If $i \neq 0$, then by (\ref{ftildei}), 
	  \begin{eqnarray}\nonumber
	  \tfi^k(b \otimes u) = (\tfi b) \otimes (\tfi^{k-1}u)
	  \mbox{ for some  }k \geq 1 
	  \end{eqnarray}
	  and $\tfi^{k-1}u \in \ovB(\lambda_{cl};0)$ hold. Together with (\ref{0apl}),
	  we may assume $b=u_{\omega_0}$.\\
	  If $\tfi u \neq 0$ for some $1 \leq i \leq N$,
	  \begin{eqnarray}\nonumber
	   \tfi(u_{\omega_0} \otimes u) = u_{\omega_0} \otimes \tfi u \neq 0.
	  \end{eqnarray}
	  This contradicts to the fact that $b_{\omega_0} \otimes u \in LW(J)$. 
	  Hence $LW(J) = \{u_{\omega_0} \otimes u_{\lambda}\}$. This
	  means that $J$ is the connected component of $B(\omega_0) \otimes B(\lambda)$ containing
         $u_{\omega_0} \otimes u_{\lambda}$.
	 \end{proof}
	\end{prop}

        For $\nu \in W$, $\lam$ as in Lemma \ref{lemshift}, the generalized
        Littlewood-Richardson rule \cite{Nakashima} gives the
        decomposition of the tensor
        product $\ovB(\nu) \otimes \ovB(\lam_{cl})$.
	In the decomposition
        \begin{equation}\label{nulam}
	 \begin{array}{ccccc}
	  \ovB(\nu)& \otimes& \ovB(\lam_{cl}) 
	   &\xrightarrow{\sim}& \bigoplus_{j=1}^{j_{\nu,\lam}} \ovB(\mu^j(\nu,\lam_{cl})),\\
	  
	  u_{\nu}&   \otimes& u_j             &\mapsto           &  w_j
	 \end{array}
	\end{equation}
	let $w_j$ be the lowest weight vector of $\ovB(\mu^j(\nu,\lam))$,
	and $u_{\nu} \otimes u_j$ ($u_j \in \ovB(\lam_{cl})$) be the corresponding vector in the
	left hand side.

	In order to to know the weight of $w_j$ as $U_q(D(N,1))$-module, we have to
	obtain the coefficients of $\ome0$.
	\begin{defn}
	 Assume $N \geq 3$ and 
	 $u_j = \widetilde{e_{l_1}} \widetilde{e_{l_2}} \cdots \widetilde{e_{l_p}} u_{\lam_{cl}}$
	 ($1 \leq l_i \leq N$, $p \in \pos$) in (\ref{nulam}). 
	 Then we define 
	 a positive integer
	 \begin{eqnarray}\label{defa}
	  a(\mu^j(\nu,\lam_{cl})) = \sharp\left\{i ; l_i = 1 \right\}.
	 \end{eqnarray}	
	See Example \ref{ex1} for the case of $N=2$.
	\end{defn}

	The decomposition of a typical $U_q(D(N,1))$-crystal as in 
	Corrolary \ref{ome} is the following proposition.
	\begin{prop}\label{lamome0}
	 Let $\lam$ be as in (\ref{lam}), $\mu^j(\nu,\lam_{cl})$ be
	 as in (\ref{nulam}). Then we have
	 \begin{eqnarray}
	  B(\lam + \omega_0) \cong \bigoplus_{\nu \in W} 
	                           \bigoplus_{n \in \mathbb{Z}_{\geq 0}} 
			           \bigoplus_{i=1}^{i_{\nu}} 
			           \bigoplus_{j=1}^{j_{\nu,\lam}}
	   \ovB\left(\mu^j(\nu,\lam_{cl}) ; z^i(\nu) + a\left(\mu^j(\nu,\lam_{cl})\right) +2n\right).
	 \end{eqnarray}
	 
	 \begin{proof}
	  By Proposition \ref{ome} and Proposition \ref{direct},
	  \begin{align*}
	   B(\lam + \omega_0) &\cong B(\ome0) \otimes \ovB(\lambda_{cl};0)\\
	                      &\cong \bigoplus_{\nu \in W}
			             \bigoplus_{i=1}^{i_{\nu}}
				     \bigoplus_{n \in \mathbb{Z}_{\geq 0}}
			              \ovB\left(\nu ; z^i(\nu) +2n\right)
				 \otimes
				     \ovB(\lambda_{cl};0)
	  \end{align*}
	  holds. Because $\widetilde{e_1}$ makes the coefficient of $\ome0$
	  increase by $1$, we get 
	  \begin{eqnarray}\nonumber
	   B(\lam + \omega_0)  \cong \bigoplus_{\nu \in W} 
	                             \bigoplus_{i=1}^{i_{\nu}}
			             \bigoplus_{n \in \mathbb{Z}_{\geq 0}} 
			             \bigoplus_{j=1}^{j_{\nu,\lam}}
	   \ovB\left(\mu^j(\nu,\lam_{cl}) ; z^i(\nu) + a\left(\mu^j(\nu,\lam_{cl})\right) +2n\right).
	  \end{eqnarray}
	  
	 \end{proof}
	\end{prop}

\section{Tensor Products of typical representations}

%
	

        Assume $\lambda' = -\sum_{i=1}^{N}n_i'\omega_i$, $n_i'\in{\mathbb{Z}}_{\geq0}$ for
         $1 \leq i \leq N$ and $\mu^j(\nu,\lam_{cl})$ as in (\ref{nulam}). We
         also assume
        \begin{eqnarray}\label{lamdnu}
	 \begin{array}{ccccc}
	  \ovB(\lam_{cl}')& \otimes& \ovB(\mu^j(\nu,\lam_{cl})) &\xrightarrow{\sim}& 
	   \bigoplus_{k=1}^{k_{\lam',\nu,\lam}}
	   \ovB(\mu_k^j(\lam'_{cl},\nu,\lam_{cl}))\\
	  u_{\lam'}&       \otimes& x_k&                     \mapsto &  y_k
	 \end{array}
	\end{eqnarray}
	holds as in (\ref{nulam}).

	We state the main theorem.
	\begin{thm}[Main Theorem]
	 Assume
	\begin{equation}\label{8.1}
	 \lam = -\sum_{i=1}^{N}n_i\omg_i,\quad
	 \lam' = -\sum_{i=1}^{N}n_i'\omg_i,\quad n_i,n_i' \in \pos \mbox{ for }1 \leq i \leq N.
	\end{equation}
	 Then we obtain
	 \begin{equation}
	  \label{main}
	   \begin{split}
	    &B(\lambda' + \omega_0) \otimes B(\lambda + \omega_0)\\
	    =& \bigoplus_{\nu \in W}
	       \bigoplus_{n \in \pos} 
	       \bigoplus_{i=1}^{i_{\nu}} 
	       \bigoplus_{j=1}^{j_{\nu,\lam}}
	       \bigoplus_{k=1}^{k_{\lam',\nu,\lam}}\\
	    &B\left(
	           \mu^j_k(\lambda_{cl}', \nu, \lambda_{cl})_{su} 
	           +\left\{z^i(\nu)+
	                   a(\mu^j(\nu, \lambda_{cl})) + 
	                   a(\mu^j_k(\lambda'_{cl}, \nu, \lambda_{cl})) + 2n + 1
	            \right\}
	            \omega_0 
	     \right).
	  \end{split}
	 \end{equation}
	 \begin{proof}
	  By Lemma \ref{lwv}, it suffices to consider
	  $u_{\lam' + \ome0} \otimes u$, $u \in B(\lam + \ome0)$.

	  Assume
	  $u \in \ovB\left(\mu^j(\nu,\lam_{cl}) ; z^i(\nu) + a\left(\mu^j(\nu,\lam_{cl})\right) +2n\right)$.
	  We have
	  \begin{equation}\nonumber
	   \begin{split}
	   &u_{\lam' + \ome0} \otimes u \in LW(B(\lambda' + \omega_0) \otimes B(\lambda + \omega_0))\\
	    \Longleftrightarrow
	   &u_{\lam' + \ome0} \otimes u \in \overline{LW}\left\{ \ovB(\lam'_{cl};1) \otimes 
            \ovB\left(\mu^j(\nu,\lam_{cl}) ; z^i(\nu) + a\left(\mu^j(\nu,\lam_{cl})\right) +2n\right)\right\}\\
	    \Longleftrightarrow
	   &u_{\lam' + \ome0} \otimes u \in \overline{LW}\left\{\ovB\left(\mu^j_k(\lambda_{cl}', \nu, \lambda_{cl});
	    z^i(\nu)+a(\mu^j(\nu, \lambda_{cl})) + a(\mu^i_k(\lam'_{cl}, \nu, \lam_{cl})) + 2n +1\right)\right\} \\
	   &\mbox{ for some }k.\\
	  \end{split}
	  \end{equation}
	  Here, the first equality follows from Lemma \ref{typeasy} and
	  Proposition \ref{lamome0}, 
	  the second
	  from (\ref{lamdnu}). Because
	  $wt(\overline{LW}(\ovB(\Lam;p\ome0)))=\Lam_{su} + p\ome0$ 
	  as $U_q(D(N,1))$-crystal, 
	  we have (\ref{main}).
	\end{proof}
	\end{thm}

	\begin{ex}\label{ex1}: $U_q(D(2,1))$

	The even part is 
	$\left(U_q(\mathfrak{sl}_2) \otimes U_q(\mathfrak{sl}_2)\right)\otimes U_q(C(1))$. 
	Let $\lam$ and $\lam'$ be as in (\ref{8.1}) and
	let $\Lam_1$ and $\Lam_2$ be the fundamental
	weights of first two $U_q(\mathfrak{sl}_2)$. 
	Because not only
	the first but also
	the second node are connected with $0$-th node in the Dynkin
	diagram, we modify (\ref{defa}) into
	\begin{equation}\nonumber
	  a(\mu^j(\nu,\lam_{cl})) = \sharp\left\{i ; l_i = 1 \mbox{ or }2\right\}.	 
	\end{equation}
	Note that
	\begin{equation}\nonumber
	 \ovB(\lam_{cl}') \otimes \ovB(\lam_{cl})
	 = \bigoplus_{j=0}^{\min(n_1,n_1')}
	   \bigoplus_{k=0}^{\min(n_2,n_2')}
	  \ovB\left(
	             -(n_1 + n_1' - 2j)\Lam_1 -(n_2 + n_2' - 2k)\Lam_2
	       \right)
	\end{equation}
	by the Clebsch-Gordan formula, and
	\begin{equation}\nonumber
	 a\left(
	        -(n_1 + n_1' - 2j)\Lam_1 -(n_2 + n_2' - 2k)\Lam_2
	  \right)
	 = j + k.
	\end{equation}
	We assume $|n_1 - n_1'| \geq 2$ and $|n_2 - n_2'| \geq 2$
	to make the description simple.

	Applying the Clebsch-Gordan formula again, we obtain
	\begin{align*}
	  &B(\lam' + \ome0) \otimes B(\lam + \ome0)\\
	 =&\bigoplus_{j=0}^{\min(n_1,n_1')}
	   \bigoplus_{k=0}^{\min(n_2,n_2')}
	   \bigoplus_{n \in \pos}\\
	  &B\left(
	          -(n_1 + n_1' - 2j)\omg_1 -(n_2 + n_2' - 2k)\omg_2
	          +(j+k+2+2n)\ome0
	   \right)\\
	 \oplus
	 &B\left(
	          -(n_1 + n_1' - 2j)\omg_1 -(n_2 + n_2' - 2k)\omg_2
	          +(j+k+6+2n)\ome0
	   \right)\\
	 \oplus
	  &B\left(
	          -(n_1 + n_1' - 2j +1)\omg_1 -(n_2 + n_2' - 2k +1)\omg_2
	          +(j+k+2+2n)\ome0
	   \right)\\
	 \oplus
	 &B\left(
	          -(n_1 + n_1' - 2j +1)\omg_1 -(n_2 + n_2' - 2k +1)\omg_2
	          +(j+k+4+2n)\ome0
	   \right)\\
	 \oplus
	 &B\left(
	          -(n_1 + n_1' - 2j +1)\omg_1 -(n_2 + n_2' - 2k -1)\omg_2
	          +(j+k+3+2n)\ome0
	   \right)\\
	 \oplus
	 &B\left(
	          -(n_1 + n_1' - 2j +1)\omg_1 -(n_2 + n_2' - 2k -1)\omg_2
	          +(j+k+5+2n)\ome0
	   \right)\\
	 \oplus
	 &B\left(
	          -(n_1 + n_1' - 2j -1)\omg_1 -(n_2 + n_2' - 2k +1)\omg_2
	          +(j+k+3+2n)\ome0
	   \right)\\
	 \oplus
	 &B\left(
	          -(n_1 + n_1' - 2j -1)\omg_1 -(n_2 + n_2' - 2k +1)\omg_2
	          +(j+k+5+2n)\ome0
	   \right)\\
	 \oplus
	 &B\left(
	          -(n_1 + n_1' - 2j -1)\omg_1 -(n_2 + n_2' - 2k -1)\omg_2
	          +(j+k+4+2n)\ome0
	   \right)\\
	 \oplus
	 &B\left(
	          -(n_1 + n_1' - 2j -1)\omg_1 -(n_2 + n_2' - 2k -1)\omg_2
	          +(j+k+6+2n)\ome0
	   \right)\\
	 \oplus
	 &B\left(
	          -(n_1 + n_1' - 2j +2)\omg_1 -(n_2 + n_2' - 2k)\omg_2
	          +(j+k+3+2n)\ome0
	   \right)\\
	 \oplus
	 &B\left(
	          -(n_1 + n_1' - 2j +0)\omg_1 -(n_2 + n_2' - 2k)\omg_2
	          +(j+k+4+2n)\ome0
	   \right)\\
	 	 \oplus
	 &B\left(
	          -(n_1 + n_1' - 2j -2)\omg_1 -(n_2 + n_2' - 2k)\omg_2
	          +(j+k+5+2n)\ome0
	   \right)\\
	 \oplus
	  &B\left(
	          -(n_1 + n_1' - 2j)\omg_1 -(n_2 + n_2' - 2k +2)\omg_2
	          +(j+k+3+2n)\ome0
	   \right)\\
	 \oplus
	 &B\left(
	          -(n_1 + n_1' - 2j)\omg_1 -(n_2 + n_2' - 2k +0)\omg_2
	          +(j+k+4+2n)\ome0
	   \right)\\
	 \oplus
	 &B\left(
	          -(n_1 + n_1' - 2j)\omg_1 -(n_2 + n_2' - 2k -2)\omg_2
	          +(j+k+5+2n)\ome0
	   \right).
	\end{align*}
	\end{ex}

	\begin{ex}: $B(-\omg_4 + \ome0) \otimes B(\ome0)$ for $U_q(D(4,1))$
	
	$W$ for $U_q(D(4,1))$ is
	\begin{equation}\nonumber
	 W = \{ \xi_0 = 0, -\xi_1 = -\Lam_1, -\xi_2 = -\Lam_2, -\xi_3 = -\Lam_3-\Lam_4,
	       -\xi_4' = -2\Lam_3, -\xi_4=-2\Lam_4 \}.
	\end{equation}
	In $U_q(D(4))$, we have
	\begin{equation}\nonumber
	 \begin{array}{lllll}
	 \ovB(-\Lam_4) \otimes \ovB(-\xi_1) 
	 =& \ovB(-\Lam_1-\Lam_4)     &\oplus \ovB(-\Lam_3),&& \\
	      &a(-\Lam_1-\Lam_4)=0,   &a(-\Lam_3)=1,&&\\
	  &&&&\\

	 \ovB(-\Lam_4) \otimes \ovB(-\xi_2) 
	 =& \ovB(-\Lam_2-\Lam_4)     &\oplus \ovB(-\Lam_1-\Lam_3) &\oplus \ovB(-\Lam_4), &\\
	      &a(-\Lam_2-\Lam_4)=0,  &a(-\Lam_1-\Lam_3)=0,        &a(-\Lam_3)=1,&\\
	  &&&&\\

	 \ovB(-\Lam_4) \otimes \ovB(-\xi_3) 
	 =& \ovB(-\Lam_3-2\Lam_4)    &\oplus \ovB(-\Lam_2-\Lam_3) &\oplus \ovB(-\Lam_1-\Lam_4) \oplus \ovB(-\Lam_3), \\
	      &a(-\Lam_3-2\Lam_4)=0, &a(-\Lam_1-\Lam_3)=0,        &a(-\Lam_1-\Lam_4)=0,  a(-\Lam_3)=1,\\
	  &&&&\\

	 \ovB(-\Lam_4) \otimes \ovB(-\xi_4') 
	 =& \ovB(-2\Lam_3-\Lam_4)    &\oplus \ovB(-\Lam_1-\Lam_3), &&\\
	      &a(-2\Lam_3-\Lam_4)=0, &a(-\Lam_1-\Lam_3)=0,&&\\
	  &&&&\\

	 \ovB(-\Lam_4) \otimes \ovB(-\xi_4) 
	 =& \ovB(-3\Lam_4)           &\oplus \ovB(-\Lam_2-\Lam_4) &\oplus \ovB(-\Lam_4), &\\
	      &a(-3\Lam_4)=0,        & a(-\Lam_2-\Lam_4)=0,       &  a(-\Lam_4)=1.&\\
	 \end{array}
	\end{equation}

	It follows that
	\begin{align*}
	 &B(-\omg_4 + \ome0) \otimes B(\ome0)\\
	 =&\bigoplus_{n \in \pos}
	  B\left(
	          -\omg_4
	          +(2+2n)\ome0
	   \right)
	 \oplus
	 B\left(
	          -\omg_4
	          +(10+2n)\ome0
	   \right)\\
	 \oplus
	  &B\left(
	          -\omg_1-\omg_4
	          +(2+2n)\ome0
	   \right)
	 \oplus
	  B\left(
	          -\omg_1-\omg_4
	          +(8+2n)\ome0
	   \right)\\
	 \oplus
	  &B\left(
	          -\omg_3
	          +(3+2n)\ome0
	   \right)
	 \oplus
	  B\left(
	          -\omg_3
	          +(9+2n)\ome0
	   \right)\\
	 \oplus
	  &B\left(
	          -\omg_2-\omg_4
	          +(3+2n)\ome0
	   \right)
	 \oplus
	  B\left(
	          -\omg_2-\omg_4
	          +(7+2n)\ome0
	   \right)\\
	 \oplus
	  &B\left(
	          -\omg_1-\omg_3
	          +(3+2n)\ome0
	   \right)
	 \oplus
	  B\left(
	          -\omg_1-\omg_3
	          +(7+2n)\ome0
	   \right)\\
	 \oplus
	  &B\left(
	          -\omg_4
	          +(4+2n)\ome0
	   \right)
	 \oplus
	  B\left(
	          -\omg_4
	          +(8+2n)\ome0
	   \right)\\
	 \oplus
	  &B\left(
	          -\omg_3-2\omg_4
	          +(4+2n)\ome0
	   \right)
	 \oplus
	  B\left(
	          -\omg_3-2\omg_4
	          +(6+2n)\ome0
	   \right)\\
	 \oplus
	  &B\left(
	          -\omg_2-\omg_3
	          +(4+2n)\ome0
	   \right)
	 \oplus
	  B\left(
	          -\omg_2-\omg_3
	          +(6+2n)\ome0
	   \right)\\
	 \oplus
	  &B\left(
	          -\omg_1-\omg_4
	          +(4+2n)\ome0
	   \right)
	 \oplus
	  B\left(
	          -\omg_1-\omg_4
	          +(6+2n)\ome0
	   \right)\\
	 \oplus
	  &B\left(
	          -\omg_3
	          +(5+2n)\ome0
	   \right)
	 \oplus
	  B\left(
	          -\omg_3
	          +(7+2n)\ome0
	   \right)\\
	 \oplus
	  &B\left(
	          -2\omg_3-\omg_4
	          +(5+2n)\ome0
	   \right)
	 \oplus
	  B\left(
	          -\omg_1-\omg_3
	          +(5+2n)\ome0
	   \right)\\
	 \oplus
	  &B\left(
	          -3\omg_4
	          +(5+2n)\ome0
	   \right)
	 \oplus
	  B\left(
	          -\omg_2-\omg_4
	          +(5+2n)\ome0
	   \right)
	 \oplus
	  B\left(
	          -\omg_4
	          +(6+2n)\ome0
	   \right).
 	\end{align*}
	\end{ex}
\section{Results for {\boldmath$U_q(B(N,1))$}}


        The results in previous sections carry over to the 
	case of $U_q(B(N,1))$.
	Since the proofs are entirely similar,
	we only state the results.

	The Cartan matrix for $U_q(B(N,1))$ is 
	\begin{eqnarray}\label{bcar}
	  A=(\left<h_i,\alpha_i\right>)_{i,j=0}^{N}=(a_{ij})_{i,j=0}^{N}= \left(\begin{array}{rrrrrrrr}

	    0       & 1      & 0      & \cdots & \\

	    -1      & 2      & -1     &\\
	    0       & -1     & 2      & -1     &\\
	    \vdots  &        & \ddots & \ddots & \ddots&       & \vdots & \\
                    &        &        & -1     & 2     & -1    & 0      & \vdots\\
                    &        &        &        & -1    & 2     & -1     &0\\
	            &        &        & \cdots & 0     & -1    & 2      &-1\\
                    &        &        &        & \cdots& 0     & -2     &2\\
	                       \end{array}\right).
	\end{eqnarray}

	The associated Dynkin diagram is \\

	\begin{center}
\unitlength 0.1in
\begin{picture}( 44.0000,  2.3000)(  5.0000, -7.3000)
%
\special{pn 8}%
\special{ar 600 600 100 100  0.0000000 6.2831853}%
%
\special{pn 8}%
\special{ar 1200 600 100 100  0.0000000 6.2831853}%
%
\special{pn 8}%
\special{ar 1800 600 100 100  0.0000000 6.2831853}%
%
\special{pn 8}%
\special{ar 4200 600 100 100  0.0000000 6.2831853}%
%
\special{pn 8}%
\special{ar 4800 600 100 100  0.0000000 6.2831853}%
\put(5.5000,-9.0000){\makebox(0,0)[lb]{$0$}}%
\put(11.5000,-9.0000){\makebox(0,0)[lb]{$1$}}%
\put(17.5000,-9.0000){\makebox(0,0)[lb]{$2$}}%
\put(47.0000,-9.0000){\makebox(0,0)[lb]{$N$}}%
\put(40.2000,-9.0000){\makebox(0,0)[lb]{$N-1$}}%
%
\special{pn 8}%
\special{sh 1}%
\special{ar 2400 600 10 10 0  6.28318530717959E+0000}%
\special{sh 1}%
\special{ar 2600 600 10 10 0  6.28318530717959E+0000}%
\special{sh 1}%
\special{ar 2800 600 10 10 0  6.28318530717959E+0000}%
\special{sh 1}%
\special{ar 3000 600 10 10 0  6.28318530717959E+0000}%
\special{sh 1}%
\special{ar 3200 600 10 10 0  6.28318530717959E+0000}%
\special{sh 1}%
\special{ar 3400 600 10 10 0  6.28318530717959E+0000}%
\special{sh 1}%
\special{ar 3600 600 10 10 0  6.28318530717959E+0000}%
\special{sh 1}%
\special{ar 3600 600 10 10 0  6.28318530717959E+0000}%
%
\special{pn 8}%
\special{pa 530 530}%
\special{pa 670 670}%
\special{fp}%
\special{pa 670 530}%
\special{pa 530 670}%
\special{fp}%
%
\special{pn 8}%
\special{pa 700 600}%
\special{pa 700 600}%
\special{fp}%
\special{pa 1100 600}%
\special{pa 700 600}%
\special{fp}%
\special{pa 1300 600}%
\special{pa 1700 600}%
\special{fp}%
\special{pa 1900 600}%
\special{pa 2300 600}%
\special{fp}%
\special{pa 3700 600}%
\special{pa 4100 600}%
\special{fp}%
\special{pa 4700 600}%
\special{pa 4600 500}%
\special{fp}%
\special{pa 4700 600}%
\special{pa 4600 700}%
\special{fp}%
\special{pa 4650 550}%
\special{pa 4280 550}%
\special{fp}%
\special{pa 4280 650}%
\special{pa 4650 650}%
\special{fp}%
\end{picture}%
 \mbox{ . }\\
	\end{center}

	We put
	\begin{equation}\nonumber
	 l_0 = 2,\quad l_1 = \cdots = l_{N-1} = -2 ,\quad l_N=-1.
	\end{equation}
	
	$U_q(B(N,1))$ is defined as in Definition \ref{yama}, wherein
	$(a_{ij})$ is replaced by
	the Cartan matrix (\ref{bcar}).

	$\{ \delta, \varepsilon_1,\ldots,\varepsilon_N \}$,
	$\Delta_0, \Delta_1, \overline{\Delta_1}$ for $U_q(B(N,1))$ is
	given by
	\begin{equation}\nonumber
	 \alpha_0=\delta-\varepsilon_1,\quad \alpha_1=\varepsilon_1-\varepsilon_2,\ldots,
	 \alpha_{N-1}=\varepsilon_{N-1}-\varepsilon_N, \quad \alpha_N=\varepsilon_N,
	\end{equation}
	\begin{equation}\nonumber
	 \Delta_0=\{\pm \varepsilon_i \pm \varepsilon_j, \pm\varepsilon_i, \pm2\delta\}, \quad 
	 \Delta_1=\{\pm \varepsilon_i \pm \delta, \pm \delta\},\quad
	 \overline{\Delta_1} = \{\pm \varepsilon_i \pm \delta\}.
	\end{equation}

	 Its even part $U_q(B(N,1)_0)$ is given by
	 $U_q(B(N)) \otimes U_q(C(1))$,
	 where
	 $U_q(B(N))$ is the subalgebra with generators 
	 $e_i, f_i,q^{h_i}$ $(1 \leq i \leq N)$,
	 and $U_q(C(1)) \simeq U_q(\mathfrak{sl}_2)$ is the one generated by
	 $E,F,q^H$, where $H=2(h_0-h_1-\cdots-h_{N-1})-h_N$,
	 $E$ and $F$ are the elements corresponding to the root
	 $2(\alpha_0 + \cdots + \alpha_{N-2}) + \alpha_N$.

%
	
	We state properties of $U_q(B(N))$ corresponding to those of
	$U_q(D(N))$ in Section 4.
	$U_q(B(N))$ has one spin representation $\ovB_{sp} =\ovB(-\Lam_N)$ 
	whose crystal base is realized as
	\begin{equation}\nonumber
	 \ovB_{sp} = \{b=(i_1,\ldots,i_N) ; i_1,\ldots,i_N=\pm\}
	\end{equation}
	with the lowest weight vector $(+,\ldots,+)$.
	The actions of $\tei$ and $\tfi$ are
	\begin{eqnarray}\nonumber
	    \widetilde{f_l}(i_1,i_2,\ldots,i_N) &=&
		    \begin{cases}
		     (i_1,\ldots,\stackrel{l}{+},\stackrel{l+1}{-},\ldots,i_N) &\mbox{ if } i_l=-,i_{l+1}=+,\\\nonumber 
		     0                       &\mbox{ otherwise, }
		    \end{cases}\\\nonumber 
	    \widetilde{e_l}(i_1,i_2,\ldots,i_N) &=&
		    \begin{cases}
		     (i_1,\ldots,\stackrel{l}{-},\stackrel{l+1}{+},\ldots,i_N) &\mbox{ if } i_l=+,i_{l+1}=-,\\\nonumber 
		     0                       &\mbox{ otherwise, }
		    \end{cases}\\\nonumber
	\end{eqnarray}
	for $1 \leq l \leq N-1$, and
	\begin{eqnarray}\nonumber
	    \widetilde{f_N}(i_1,i_2,\ldots,i_N) &=&
		    \begin{cases}
		     (i_1,\ldots,\stackrel{N}{+}) &\mbox{ if } i_N=-,\\\nonumber 
		     0                      &\mbox{ otherwise, }
		    \end{cases}\\\nonumber 
	    \widetilde{e_N}(i_1,i_2,\ldots,i_N) &=&
		    \begin{cases}\nonumber 
		     (i_1,\ldots,\stackrel{N}{-}) &\mbox{ if } i_N=+,\\\nonumber 
		     0                      &\mbox{ otherwise. }
		    \end{cases}
	 \end{eqnarray}
	 
	 We also describe $\ovB_{sp}$ in terms of Young tableaux as we do
	 in the case of $U_q(D(N))$.
	We denote
	\begin{equation}\nonumber
	 \xi_0=0, \quad \xi_i = \Lam_i \quad(1 \leq i \leq N-1),\quad \xi_N=2\Lam_N.
	\end{equation}

	\begin{prop}
        We have
	 \begin{eqnarray}\nonumber
	  \ovB_{sp} \otimes \ovB_{sp} = \bigoplus_{0 \leq k \leq N}
	   \ovB(-\xi_k).
	 \end{eqnarray}
	 For $0 \leq i \leq N$, the lowest weight vector
	 corresponding to the connected component $\ovB(-\xi_i)$ is 
	 \begin{eqnarray}\nonumber 
	  (+,\ldots,+) \otimes (\overbrace{+,\ldots,+}^{i},-,\ldots,-).	   
	 \end{eqnarray}	 
	\end{prop}

	\begin{prop}[Koga\cite{Koga}]
	 \begin{description}
	  \item[(1)]Assume  $u \otimes v = t(a_1,\ldots,a_N) \otimes t(b_1,\ldots,b_N)\in \ovB_{sp} \otimes \ovB_{sp}$. 
	 Then we have\\
	 \begin{description}
	  \item[(1A)]For $0 \leq k \leq N-1$,
	 \begin{align*}  
	     &u \otimes v \in \ovB(-\xi_k)  
	     \Longleftrightarrow\\
	     &t(a_1,\ldots,a_{N-k};a_{N-k+1},\ldots,a_N|b_1,\ldots,b_k;b_{k+1},\ldots,b_N)
               \mbox{ is semi-standard and }\\
	     &t(a_1,\ldots,a_{N-k-1};a_{N-k},\ldots,a_N|b_1,\ldots,b_{k+1};b_{k+2},\ldots,b_N)
	  \mbox{ is not semi-standard,}	 
	 \end{align*}
	  \item[(1B)]
	 \begin{eqnarray}\nonumber  
	     u \otimes v \in \ovB(-\xi_N)  
	     \Longleftrightarrow
	     t(;a_1,\ldots,a_N|b_1,\ldots,b_N;)
	     \mbox{ is semi-standard. }
	 \end{eqnarray}
	 \end{description}
	\end{description}
	\end{prop}

	$U_q(B(N,1))$-module $V(-\omg_N)$ admits a crystal base.

	\begin{prop}
		The irreducible lowest weight module $V(-\omg_N)$ with lowest weight
	        $-\omega_N$ has a basis over ${\mathbb{Q}}(q)$
	 
		\begin{equation}\nonumber
		\{{v}(i_1,\ldots,i_N)_k \>; \> k \in \pos, (i_1,\ldots,i_N)
		\in \ovB_{sp}\}
		\end{equation}
		with the lowest weight vector ${v}(+,\ldots,+)_0$
	        such that the actions of $\sigma$ and $e_i$ are;
	        \begin{eqnarray}\nonumber
		 \sigma {v}(+,\ldots,+)_0 = {v}(+,\ldots,+)_0, 
		\end{eqnarray}
		\begin{eqnarray}\nonumber
		  {e_i}({v}(i_1,\ldots,i_N)_k) &=&
		    \begin{cases}
		     {v}(i_1',\ldots,i_N')_k &\mbox{ if }\widetilde{e_i}(i_1,\ldots,i_N)=
		                                              (i_1',\ldots,i_N') \neq 0 \text{ in }\ovB_{sp},\\\nonumber 
		     0                             &\mbox{ otherwise, }
		    \end{cases}\\\nonumber
		\end{eqnarray}
		for $1 \leq i \leq N$, and
		\begin{eqnarray}\nonumber
		  {e_0}({v}(i_1,i_2,\ldots,i_N)_k) &=&
		    \begin{cases}
		     q^{-k}v(+,i_2,\ldots,i_N)_{k+1} &\mbox{ if }i_1=-,\\\nonumber 
		     0                                 &\mbox{ otherwise. }
		    \end{cases}
		\end{eqnarray}
	 \end{prop}

	\begin{prop}
		The irreducible lowest weight module $V(-{\omega}_N)$ has the polarizable crystal base $(L,B)$;
		\begin{eqnarray}\nonumber
		L = 	\bigoplus_{\substack{(i_1,\ldots,i_N) \in \ovB_{sp}\\
				k \in \pos }}
				Av(i_1,\ldots,i_N)_k,
		\end{eqnarray}
		\begin{eqnarray}\nonumber
		B =& \{\pm v(i_1,\ldots,i_N)_k \mod qL ; (i_1,\ldots,i_N) \in \ovB_{sp}, k \in \pos\}.
		\end{eqnarray}
	 
	 The Kashiwara operators on $B$ is given by (we omit $\mod qL$)
	 \begin{eqnarray}\nonumber
	    \widetilde{e_i}v(i_1,\ldots,i_N)_k&=&
		    \begin{cases}
		     v(i_1',\ldots,i_N')_k &\mbox{ if }\widetilde{e_i}(i_1,\ldots,i_N)=
		                                              (i_1',\ldots,i_N') \neq 0 \text{ in }\ovB_{sp},\\\nonumber 
		     0                                  &\mbox{ otherwise, }
		    \end{cases}\\\nonumber
	 \end{eqnarray}
		for $1 \leq i \leq N$, and 
	 \begin{eqnarray}\nonumber
	    \widetilde{e_0}v(i_1,i_2,\ldots,i_N)_k&=&
		    \begin{cases}
		     v(+,i_2,\ldots,i_N)_{k+1} &\mbox{ if }i_1=-,\\\nonumber 
		     0                                  &\mbox{ otherwise. }
		    \end{cases}
	 \end{eqnarray}
	\end{prop}

	Because
	 \begin{eqnarray}\nonumber
	  B(-\omega_N) \otimes B(-\omega_N) = B(-2\omega_N) \oplus \bigoplus\limits_{j=1}^{N-1}B(-\omega_j)
	   \oplus \bigoplus\limits_{k \in \pos}B((k+1)\omega_0)
	 \end{eqnarray}
	 holds, we have the theorem corresponding to Theorem \ref{exist}.
	\begin{thm}
	 The irreducible lowest weight module $V(\lambda)$ with
	 the lowest weight
	 \begin{equation}\nonumber
		\lambda = n_0{\omega}_0 - \sum_{i=1}^{N}n_i{\omega}_i, \quad n_i \in \pos \quad \mbox{ for }
		 0 \leq i \leq N
	 \end{equation} 
        admits a polarizable crystal base.
	\end{thm}

	The $0$-arrows in $B(-\omega_N) \otimes B(-\omega_N)$ can be
	described as follows.

	\begin{lem}\label{0connection}
	 In $B(-\omega_N) \otimes B(-\omega_N)$, we have\\
	 \begin{description}
	 \item[(1)]for $k=2,3,\ldots,N$
	 $$\begin{CD}
	  \ovB(-\xi_k;l) @>{0}>R> \ovB(-\xi_{k-1};l-1),
	 \end{CD}$$

	 \item[(2)]
	 $$\begin{CD}
	    \hspace*{-8mm}
	  \ovB(-\xi_1;l) @>{0}>R> \ovB(-\xi_0;l),
	 \end{CD}$$

	 \item[(3)]for $k=2,\ldots,N$
	$$\begin{CD}
	   \hspace*{-5mm}
	  \ovB(-\xi_{k-1};l) @>{0}>L> \ovB(-\xi_k;l-1),
	 \end{CD}$$
	
	 \item[(4)]
	$$\begin{CD}
	  \ovB(-\xi_0;l) @>{0}>L> \ovB(-\xi_1;l-2),
	 \end{CD}$$

	 \item[(5)]
	 $$\begin{CD}
	  \ovB(-\xi_N;l) @>{0}>L> \ovB(-\xi_N;l-1).
	   \end{CD}$$
	 \end{description}
	\end{lem}

	The results corresponding to Corollary \ref{ome} is the next proposition.
	\begin{prop}
	 We have
	 \begin{eqnarray}\nonumber
	  B(\omega_0)=\bigoplus_{\nu \in W}\bigoplus_{i=1}^{i_{\nu}}\bigoplus_{n \in \mathbb{Z}_{\geq 0}}
	   \ovB\left(\nu ; z^i(\nu) +2n\right),
	 \end{eqnarray}
	 where
	 $W = \left\{-\xi_0,\ldots,-\xi_N\right\}$. 
	 $z^i(\nu)$ are given as follows;
	 
	  \begin{equation}\nonumber
	   \begin{array}{lllll}
	    \mbox{ for } \nu = -\xi_0,&  &i_{\nu} = 2,&  z^1(\nu) = 1 , &z^2(\nu) = 2N+2,\\
	    \mbox{ for } \nu = -\xi_k \quad (1 \leq k \leq N),& 
	                       &i_{\nu} = 2,& z^1(\nu) = k , &z^2(\nu) = 2N-k+1.
	   \end{array}
	  \end{equation}
	 \end{prop}

	 The results in Section 7 and Section 8 hold for $U_q(B(N,1))$.


\vspace{7mm}
\leftline{\large\textbf{Acknowledgment}}

        The author would like to express gratitude to Professor
	M. Jimbo for his continuous guidance and support during the course of
	this work. 
	He is also grateful to S. Ariki, M. Kashiwara, 
	K. Koga and T. Nakashima for interest in this work and 
	a number of valuable comments.
	He is also indebted to H. Kajiura, S. Kato and Y. Terashima 
	for help with references.



\begin{thebibliography}{99}



\bibitem{BK}
	G.~Benkart and S.-J.~Kang,
	{\it Crystal bases for quantum superalgebras},
	Adv. Stud. Pure Math., {\bf 28}(2000), 21-54.

\bibitem{BKD}
	G.~Benkart, S.-J.~Kang,~and D.~Melville,
	{\it Quantized enveloping algebras for Borcherds superalgebras},
	J. of Amer. Math. Soc., {\bf 350(8)}(1998), 3297-3319.

\bibitem{BKK}
	G.~Benkart, S.-J.~Kang,~and M.~Kashiwara,
	{\it Crystal bases for the quantum superalgebra $U_q(\mathfrak{gl}(m,n))$},
	J. of Amer. Math. Soc., {\bf 13(2)}(2000), 295-331.

\bibitem{J}
	K.~Jeong,
	{\it Crystal bases for Kac-Moody Superalgebras},
	J. Algebra, {\bf 237}(2001), no. 2, 562-590.

\bibitem{Kac1}
	V.~Kac,
	{\it Lie superalgebras},
	Adv. Math., {\bf 26}(1977), 8-96.

\bibitem{Kac2}
	V.~Kac,
	{\it Representations of classical Lie superalgebras},
	Lecture Notes in Math., vol. 676, Springer, Berlin, 1978, 597-626.

\bibitem{KN}
	M.~Kashiwara,~and T.~Nakashima,
	{\it Crystal graphs for representations of the $q$-analogue of classical Lie algebras},
	J. Algebra, {\bf 165}(1994), 295-345.

\bibitem{Koga}
	Y.~Koga,
	{\it Level one perfect crystals for $B_n^{(1)}$,$C_n^{(1)}$,$D_n^{(1)}$},
	J. Algebra, {\bf 217}(1999), 312-334.

\bibitem{MZ}
	I.~M.~Musson and Y.-M.~Zou,
	{\it Crystal bases for $U_q({\frak o}{\frak s}{\frak p}(1,2r))$},
	J. Algebra, {\bf 210}(1998), 514-534.

\bibitem{Nakashima}
	T.~Nakashima,
	{\it Crystal base and a generalization of the Littlewood-Richardson rule for the classical Lie algebras},
	Comm. Math. Phys., {\bf 154}(1993), 215-243.

\bibitem{Yamane}
	H.~Yamane,
	{\it Quantized enveloping algebras associated with simple Lie 
	superalgebras and their universal $R$-matrices},
	Publ. RIMS., {\bf 30}(1994), 15-87.

\bibitem{Z1}
	Y.~M.~Zou, 
	{\it Crystal bases for $U_q(\Gamma(\sigma_1,\sigma_2,\sigma_3))$}, 
	Trans. Amer. Math., {\bf 353(9)}(2001), 3789-3802.

	



	




\end{thebibliography}
\end{document}